\magnification=1200

\hsize 17truecm
\vsize 22truecm

\font\ten=cmbx10

\font\twelvec=msbm10 at 12pt
\font\sevenc=msbm10 at 9pt
\font\fivec=msbm10 at 7pt

\newfam\co
\textfont\co=\twelvec
\scriptfont\co=\sevenc
\scriptscriptfont\co=\fivec

\def\out{\mathop{\rm out}\nolimits}
\def\i{\mathop{\rm in}\nolimits}
\def\im{\mathop{\rm Im}\nolimits}
\def\re{\mathop{\rm Re}\nolimits}

\def\id{\mathop{\rm Id}\nolimits}
\def\lim{\mathop{\rm lim}\nolimits}

\def\mod{\mathop{\rm mod}\nolimits}
\def\supp{\mathop{\rm supp}\nolimits}

\def\inf{\mathop{\rm inf}\nolimits}

\def\fed{\mathop{\rm def}\nolimits}
\def\diag{\mathop{\rm diag}\nolimits}

\def\neigh{\mathop{\rm neigh}\nolimits}

\def\Sum{\displaystyle\sum}
\def\e{\mathop{\rm \varepsilon}\nolimits}

\baselineskip 15pt
%\end

\centerline{\ten INTEGRABILITY, HYPERBOLIC FLOWS} 
\centerline{\ten AND THE BIRKHOFF NORMAL FORM}

\bigskip

\centerline { Michel ROULEUX }
\textfont\co=\fivec 
\centerline {Centre de Physique Th\'eorique}

\centerline {Unit\'e Propre de Recherche 7061}
\centerline {CNRS Luminy, Case 907, 13288 Marseille Cedex 9, France}

\centerline {and PhyMat, Universit\'e de Toulon et du Var}
\textfont\co=\twelvec

\bigskip

\noindent {\bf Abstract}: We prove that a Hamiltonian 
$p\in C^\infty(T^*{\bf R}^n)$ is locally integrable near a 
non-degenerate critical
point $\rho_0$ of the energy, provided that the fundamental matrix
at $\rho_0$ has no purely imaginary eigenvalues.  
This is done by using Birkhoff normal forms,
which turn out to be convergent in the $C^\infty$ sense. 
We also give versions of the Lewis-Sternberg normal form
near a hyperbolic fixed point of a canonical transformation,
using a recent result of A.Banyaga, R.de la Llave and C.Wayne.
Then we investigate the complex case, showing that when 
$p$ is holomorphic near $\rho_0\in T^*{\bf C}^n$, then 
$\re p$ becomes integrable in the complex domain for
real times, while the Birkhoff series and the Birkhoff transforms
may not converge, i.e. $p$ may not be integrable.

\bigskip

\noindent{\bf 0. Introduction}.

\medskip

Birkhoff theorem reduces hamiltonians near an elliptic 
equilibrium to quasi-integrable systems. More precisely, let 
$p \in C^\infty (T^* {\bf R}^n)$ have a local non degenerate 
minimum at $\rho_0 = (x_0, \xi_0) = 0$ with non resonant frequencies
$\lambda_1, \cdots, \lambda_n$, i.e. the fundamental matrix $F_{\rho_0}$
defined by :
$$p''_{\rho_0}(t, s)={1\over 2}\sigma(t, F_{\rho_0}(s))\leqno(0.1)$$
(here the hessian $p''$ and the symplectic 2-form are considered 
as quadratic forms on ${\bf R}^{2n}$,~) has eigenvalues $\pm i\lambda_1,
\cdots, \pm i\lambda_1$ linearly independent over ${\bf Z}$, $\lambda_j>0$,
then there is
(locally near $\rho_0$,) a canonical transform $\kappa \in C^\infty$
preserving the origin $\rho_0=0$ formally defined through its Taylor 
series, such that 
$$q(y, \eta) = p \circ \kappa (y, \eta) \sim
\Sum_{\alpha \in {\bf N}^n \setminus 0} a_\alpha \iota^\alpha ,
\ \iota_j = {1 \over 2}(\eta_j^2 + y_j^2) \leqno(0.2)
$$ 
near 0 (in the sense of Taylor series, )  with linear part
$\Sum_{j=1}^n \lambda_j \iota_j$. Function $q$ is known as the
Birkhoff normal form of $p$ (see [Bi], [Ga], [GiDeFoSi], [Sj4], [Vi], etc...) 
A theorem of C.Siegel [Si] says that
Birkhoff series are in general divergent (because of small
denominators) and there is no hope to reduce $p$ to a completely
integrable system. A gigantic litterature
has been devoted to integrability of hamiltonian systems ; we have 
listed below some of the most famous references ([Ar], [ArNo], 
[CuB], [Ga], [Mo], 
[Si], [SiMo], \dots ) but this work has been in part inspired by [El],
and [It]. See also [Au] for a somewhat less
conventionnal and more algebraic approach.

Classification of quadratic
hamiltonians was made by Williamson cf. [Ar.App.6]. We know that
eigenvalues of $F_{\rho_0}$ are of the form $\lambda, \overline \lambda,
-\lambda, -\overline \lambda$. These hamiltonians
have a
particular simple normal
form when the eigenvalues are all distinct, and non vanishing. 
Assuming that $F_{\rho_0}$ is semi-simple (diagonalizable, ) 
in suitable symplectic coordinates $(x, \xi) \in {\bf R}^{2n}$, 
the normal form 
is given as follows:
$$\eqalign{
p(x, \xi)=& \Sum_{j=1}^\ell a_j x_j\xi_j 
+\Sum_{j=1}^{m} \Bigl( c_j\bigl( x_{\ell+2j-1}
\xi_{\ell+2j-1}+x_{\ell+2j}\xi_{\ell+2j} 
\bigr)\cr
&+d_j \bigl(x_{\ell+2j-1}\xi_{\ell+2j}-x_{\ell+2j}\xi_{\ell+2j-1}\bigr)\Bigr)
+{1 \over 2}\Sum_{j=\ell+2m}^{n}b_j(\xi_j^2 + x_j^2)\cr
}\leqno(0.3)
$$
We call ''action variables'' the elementary polynomials that enter
the expression (0.3). The eigenvalues  
$\lambda_j$ of $F_{\rho_0}$ are of the form $\pm a_j$, $\pm (c_j\pm id_j)$,
and $\pm ib_j$, with the convention $a_j, \ b_j, c_j>0$. 
Here we
consider the case where none of the eigenvalues  
$\lambda_j$
is purely imaginary, i.e. no $b_j$ occur in the decomposition. 
We say then that $p$, or $H_p$ (the
hamiltonian vector field, ) is hyperbolic, or of complex hyperbolic type,
if we want to stress that some $\lambda_j$'s are complex.

Since the construction 
of Birkhoff series is a purely algebraic algorithm, it extends
trivially to the hyperbolic, or complex hyperbolic
case (provided, of course, the eigenvalues are
rationally independent.) It is commonly believed that the process
``converges'' in this situation, and the main purpose of that
paper is to provide a proof for such a  result.  

Complex eigenvalues occur in small
oscillations around an instable equilibrium. As a first example
we consider a top spinning around its apex O, 
with inertial momenta $I_1\leq I_2 < I_3$, the 
principal axis of inertia corresponding to eigenvalue $I_3$
goes through O. 
For $I_1=I_2$ (the so-called Lagrange top,~) the hamiltonian is integrable,
at all energies, but in general
there are only 2 integrals of motion. See e.g. [Au] for details.
When the top is spinning fast enough, the total energy is close
to a minimum, and the
hamiltonian orbits (expressed in suitable Euler angles) are confined 
within compact energy surfaces, on quasi-invariant torii~;
then the motion can be described
by means of the Birkhoff normal form (0.2).
Some of these torii are invariant (the KAM torii,~) but most of them
will be eventually destroyed.
When kinetic energy decreases however, 
we approach a critical value of the hamiltonian,
and the motion becomes unstable.

As a second example, we may consider a satellite, with 
inertia momenta $I_1<I_2 < I_3$, spinning around the 
principal axis of inertia corresponding to the intermediate
eigenvalue $I_2$. Again, whithin certain regimes, 
such a motion is unstable.

Then we may ask whether the hamiltonian becomes integrable near
such critical energies. 
From the point of vue of Classical Mechanics, this matter
is rather futile, since the system will leave the unstable position
long before the effects of non integrablity become relevant : 
divergence from equilibrium grows in general exponentially fast with
time, with exception however of the trajectories sufficiently
close to the stable manifold. Thus, such an improvement may be 
of ``microlocal'' nature.   

In (semi-classical) Quantum Mechanics however, 
particles are reputated to tunnel
in classically forbidden regions. A local minimum of the classical
hamiltonian becomes a saddle point 
``seen from the complex side''. 
Consider for instance 
a semiclassical Schr\"odinger operator $P = -h^2 \Delta +V(x)$ 
for energies $E$ close to a non-degenerate minimum of $V$,
$V(x_0)=0$. The classical hamiltonian reads
$p(x, \xi)= \xi^2+V(x)$.  When extending 
quasi-invariant tori in $V(x) > E$, we replace $p$ by
$\widetilde p(x, \xi)= \xi^2 - V(x)$, which becomes hyperbolic,
and it is very convenient to know, in tunneling problems 
(as in [Ro1]) that the resulting
hamiltonian, written in (hyperbolic) action-angle coordinates
is completely integrable.
Complex eigenvalues are also met when studying magnetic
Schr\"odinger operator $P(x, hD)=(hD-A(x))^2 + V(x)$ (see [MaSo], [KaRo], 
etc\dots~)

Our main result for integrability and Birkhoff transformations
in the real $C^\infty$ sense
is the following :

\smallskip
\noindent {\bf Theorem 0.1}: Assume $p \in C^\infty$ is real and
(complex-) hyperbolic with
a non-degenerate critical point at $\rho_0$, and the eigenvalues 
$\lambda_1, \cdots, \lambda_n$ are rationally
independent.  
Then there is a (germ of) $C^\infty$ canonical map $\kappa$, 
$\kappa(0, 0)=(0, 0), \ d\kappa(0, 0) = I$, and a
$C^\infty$ function $q$ of the elementary action 
variables $\iota$ as in (0.3) such that $p \circ \kappa(y, \eta) =
q(\iota)$~; the quadratic part of $q$ is as in (0.3), without elliptic
terms.
\medskip
As a by-product, we can study integrability in the neighborhood 
of a closed trajectory of hyperbolic type, as in the
examples above.

A related problem concerns conjugation of a real canonical transformation
to a time-one hamiltonian flow~; this is the so-called 
Lewis-Sternberg normal form [St]. A typical situation is 
this of the Poincar\'e map, and a lot of work has been devoted  
to the subject [Br], [Fr], [BaLlWa], [It], [IaSj] \dots.

As for the Birkhoff normal form, a central question is 
convergence of the process of reduction. The Lewis-Sternberg theorem
was stated at the level of formal series, and a proof of convergence
in the symplectic, hyperbolic case was only recently given by A. Banyaga, 
R. de la Llave and C. Wayne [BaLlWa].  

So let $\Phi : T^*{\bf R}^n \to T^*{\bf R}^n$ be a local diffeomorphism
preserving the symplectic structure, $\Phi(0, 0)=(0, 0)$. 
Assume that $d\Phi(0, 0)$ has eigenvalues $\lambda_1, \cdots, \lambda_n$,
none of them is of modulus 1. We say then that $\Phi$
is hyperbolic at (0, 0). 

Assume also the frequencies $\lambda_1, \cdots, \lambda_n$
are non resonant (in the strong sense), i.e. 
$\lambda_1^{m_1}\cdots\lambda_n^{m_n}=1$ for $m_j \in {\bf Z}$
implies $m_j=0$.

Note that if $H_p$ is a  hamiltonian
vector field, then $H_p$ is hyperbolic in the sense above,
iff the time-one map 
$\exp H_p$ is hyperbolic, because of the formula
$\kappa \circ \exp H_p \circ \kappa^{-1} = \exp H_{p\circ \kappa^{-1}}$.
Loosely speaking,
a Birkhoff normal form for $p$ gives a Sternberg normal form for 
$\exp H_p$. This is the main idea in the following~: 
\medskip
\noindent {\bf Theorem 0.2}: Let $\Phi$ be as above. Then there
is a smooth function $q(\iota)$ depending on the action variables 
$\iota$ alone,
and a smooth canonical map $\kappa$,
$\kappa(0, 0)=(0, 0)$, $d\kappa(\rho_0)=I$ such that 
$\kappa\circ\Phi\circ\kappa^{-1}(x, \xi)=\exp q(\iota)$.
\medskip
 
Next we turn to the holomorphic case, and focus on the reduction of
hamiltonians (see [It] for a discussion on necessary and sufficient
conditions ensuring that such hamiltonians are integrable. )

Again the problem arises
naturally in semi-classical Quantum Mechanics. 
As an example, consider $p(x, \xi)$ real analytic near $\rho_0=(0, 0)
\in {\bf R}^{2n}$,
with a non-degenerate minimum at $\rho_0$, and let $\pm i\lambda_1, \cdots,
\pm i\lambda_n$ be the purely imaginary eigenvalues of $F_{\rho_0}$,
$\lambda_j >0$, that  we assume again rationally independent.
When trying to construct the solution of some eikonal equation, we introduce 
$\widetilde p(z, \zeta) = -p(z-\zeta, i \zeta)$ as an holomorphic function
on a neighborhood of 0 in $T^*{\bf C}^n$.
Then $\widetilde p$
verifies the
hypotheses above, namely if $\widetilde p_2$ denotes the quadratic
part of $\widetilde p$, then
$\langle d \widetilde p_2(0, 0), (z, \zeta) \rangle
= \Sum_{j=1}^n \lambda_j z_j \zeta_j$. This situation is met when
studying microlocal properties of eigenfunctions for a magnetic
Schr\"odinger operator $P(x, hD)=(hD-A(x))^2 + V(x)$ (see [MaSo].)

As usual in complex symplectic geometry, it is convenient to
distinguish between several symplectic structures~;
we send the reader to [Sj1], [MeSj] for the theory,
and recall here simply the following fact:  ${\bf C}^{2n}$
is endowed with the complex canonical 2-form $\sigma_{\bf C} = \Sum_{j=1}^n
d\zeta_j \wedge dz_j$, $z_j=x_j+iy_j$, $\zeta_j=\xi_j+i\eta_j$,
which makes it a symplectic space, and 2 real symplectic 2-forms~: 
Re$\sigma_{\bf C}= \Sum_{j=1}^n
d\xi_j \wedge dx_j - d\eta_j \wedge dy_j$, and 
Im$\sigma_{\bf C}= \Sum_{j=1}^n
d\xi_j \wedge dy_j + d\eta_j \wedge dx_j$. 
Concerning integrability in the complex domain,
we are led  naturally to introduce the following~:
\medskip
\noindent {\bf Definition 0.3}: Let $p(x, \xi)$ be a complex 
hamiltonian near 
$\rho_0$ and have a non degenerate critical point at $\rho_0$. We say
that $p$ is R-integrable iff there is a $\re \sigma_{\bf C}$-canonical
map $\kappa \in C^\infty$ around $\rho_0$ and a $C^\infty$ function 
$q(\iota')$ such that
$\re p\circ \kappa(x, \xi)=q(\iota')$.
(Here $\iota'$ stand for the real and imaginary part of the 
complex action variables as in (0.3), and Poisson commute for the
real symplectic structure. )
\medskip
Equivalently, there exists a
$\im \sigma_{\bf C}$-canonical
map $\widetilde \kappa \in C^\infty$, and a $C^\infty$ function 
$\widehat q(\iota')$, such that
$\im p \circ \widetilde \kappa(x, \xi)=\widetilde q(\iota')$.
We could define analogously a I-integrable hamiltonian, by requiring
that $\im p \circ \kappa(x, \xi)=q(\iota')$, for some
$\re \sigma_{\bf C}$-canonical map $\kappa$. 
Roughly speaking, a R- (resp. I-) integrable hamiltonian is integrable for
real (resp. imaginary) times. If $p$ is holomorphic and ${\bf C}$-integrable, 
(i.e. with respect to $\sigma_{\bf C}$,~) then it is 
both R and I-integrable, but there are not so many
hamiltonians because of Siegel's result. In fact, H. Ito [It]
has proved that Birkhoff series and Birkhoff
transforms are convergent iff the hamiltonian is integrable
in the usual sense, e.g. the corresponding dynamical system
has, locally, $n$ Poisson commuting integrals of motion.
We have~:
\smallskip
\noindent {\bf Theorem 0.4}: Let $p(x, \xi)$ be a complex 
hamiltonian near 
$\rho_0$ and have a non degenerate critical point at $\rho_0$.
Assume that 
$\overline \partial_{(y, \eta)} p = {\cal O}(|y, \eta|^\infty)$,
and that the fundamental matrix $F_{\rho_0}$ (in the holomorphic sense)
has no purely imaginary eigenvalues.  
Then $p$ is R-integrable in a complex neighborhood of $\rho_0$. Moreover,
if $\kappa$ denotes the $\re \sigma_{\bf C}$-canonical map as in 
Definition 0.3, we have 
$\overline \partial_{(y, \eta)} \kappa = {\cal O}(|y, \eta|^\infty)$, and
$\kappa^*(\sigma_{\bf C})=\sigma_{\bf C}+{\cal O}(|y, \eta|^\infty)$.
\medskip
Our result still looks quite poor, in the sense that 
we loose on the way almost every track of analyticity~;
reduction to the normal Birkhoff form holds only modulo functions with
$\overline \partial$ of rapid decrease near $\rho_0$.
Of course again, we cannot expect convergence of Birkhoff series or Birkhoff
transforms
in a full complex neighborhood
of $\rho_0$, except in the one dimensional case, 
see [It] and [HeSj2,App.b].
A more thorough approach should rely on resurgence theory for functions
of several complex variables as in [Ec]~; 
this would of course help to understand
better how
does the system switch from integrability to non-integrability when moving 
around the origin in complex directions
(see also [Ro2] for another
type of results, where we study integrability and monodromy of $\kappa$,
as a map defined on the covering of $T^*{\bf C}^n\setminus \rho_0$, 
in the complement of the stable and unstable manifolds.)
\medskip
The paper is organized as follows:

In Section 1 we prove theorem 0.1 for hamiltonians, and give an 
equivalent formulation, via action-angle variables, which turned out 
to be useful
in computing tunneling effects for a semi-classical
Schr\"odinger operator (see [Ro1]. ) We discuss also briefly the
case of a more general center manifold.

Section 2 is devoted to the Lewis-Sternberg normal form
for canonical transforms. 

In Section 3,  we extend the Birkhoff normal form of 
theorem 0.1 
and the Sternberg normal form of theorem 0.2, to the parameter dependent case,
in the spirit of [IaSj]. 

In Section 4, we recall some wellknown facts
about complex symplectic geometry and
prove first the center
stable/unstable manifold theorem in the almost holomorphic case.
Then we turn to the proof of theorem 0.4, which is very similar
to that of theorem 0.1. We conclude with some remark on monodromy.

In the Appendix, we recall a simple way of constructing 
Birkhoff series, including parameters.

Our results have a natural extension to semi-classical
quantization as in [IaSj], but this will not be investigated here. 
We close this Introduction by listing some open problems :
\smallskip
\noindent 1) What can be said about integrability when
Spec $F_{\rho_0}\cap i{\bf R}= \{i\lambda, -i\lambda\}$,
$\lambda>0$, i.e. when the center-manifold associated with purely
imaginary eigenvalues is of dimension 2 ? For higher dimensions,
it is known that KAM torii can occur (see [Gr].)

\noindent 2) What can be said about integrability in the (complex-)
hyperbolic case, 
when
some of the frequencies are resonant, or more precisely when 
the equilibrium point $\rho_0$
is ``simply resonant'' in the sense of [It] ?

\noindent 3) Do our results extend to time-dependant hamiltonians
(see again [It]) ?
\medskip
\noindent {\bf Acknowledgements}: I want to thank J. Sj\"ostrand who gave 
inspiration
to this work, W. Craig and R. de la Llave for useful discussions.
Cheerful thanks also to I.M. Sigal for his kind hospitality at the University
of Toronto, where part of this work was done in fall 2000
under NSERC grant.

\medskip
\noindent{\bf 1. Birkhoff normal form and integrability : the
real case}.
\medskip
We discuss here of ``convergence'' of
Birkhoff normal forms for smooth hamiltonians near a fixed point $\rho_0$, 
or a closed trajectory of (possibly complex) hyperbolic type. 
\medskip
\noindent {\bf a) The hyperbolic fixed point}.
\smallskip
Let $p$ be a real valued hamiltonian with a nondegenerate critical point
$\rho_0 \in T^*{\bf R}^n$ of complex hyperbolic type. 
First we recall some wellknown facts about
the geometry of bicharacteristics of $p$ near $\rho_0$ (see [Ch], 
[Sj2], [LasSj], etc
\dots) (though 
there seems to be a confusion in [Ch,p.707], between the invariant
manifolds for the vector field $X$ and its linear part $X_0$,
the main arguments show up already in that paper.~) 
Then we discuss a solvability problem for $H_p$ in the 
class of smooth, flat functions at $\rho_0$. At last we prove Theorem 0.1.
by the method of homotopy. 
\medskip
Let
$F_{\rho_0}$ denote the
fundamental matrix of $p$ at
$\rho_0 = (0, 0)$,
$${1\over 2}F_{\rho_0}= \normalbaselineskip=18pt \pmatrix
{ {\partial^2p \over \partial x \partial \xi}&{\partial^2p \over
\partial
\xi^2}
\cr -{\partial^2p \over \partial x^2}&-{\partial^2p \over \partial x
\partial
\xi}\cr }(\rho_0) = J \ \hbox{Hess}(p)(\rho_0)
\leqno(1.1)$$
(where $J$ is the symplectic matrix, ) verifying
$\hbox{Hess}(p)(\rho_0) = p''_{\rho_0}(t, s) = {1\over 2}\sigma(t,
F_{\rho_0}(s))$. The factor ${1\over 2}$ is for convenience of notations.
Since
$p''_{\rho_0}$ is non degenerate, $F_{\rho_0}$ has no
zero eigenvalues. 
As we are interested in the Birkhoff normal form, 
we readily assume that $F_{\rho_0}$ is diagonalizable.
Let $\Lambda_\pm \subset T_{\rho_0}{\bf R}^{2n}$ be the sum of all eigenspaces
corresponding to eigenvalues with positive (resp. negative) real
parts. 

Classification of quadratic
hamiltonians was made by Williamson cf. [Ar.App.6]. We know that
eigenvalues of $F_{\rho_0}$ are of the form $\lambda, \overline \lambda,
-\lambda, -\overline \lambda$. These hamiltonians
have a
particular simple normal
form when the eigenvalues are all distinct, and non vanishing. 
Assuming that $F_{\rho_0}$ has no purely imaginary eigenvalues, 
in suitable symplectic coordinates $(x, \xi) \in {\bf R}^{2n}$, the 
normal form for the quadratic part $p_2$ of $p$ at $\rho_0$
is given as follows:
$$\eqalign{
p_2(x, \xi)&= \Sum_{j=1}^\ell a_j x_j\xi_j 
+\cr
&\Sum_{j=1}^{m} \Bigl( c_j\bigl( x_{\ell+2j-1}
\xi_{\ell+2j-1}+x_{\ell+2j}\xi_{\ell+2j} 
\bigr)
+d_j \bigl(x_{\ell+2j-1}\xi_{\ell+2j}-x_{\ell+2j}\xi_{\ell+2j-1}\bigr)\Bigr)
\cr}
\leqno(1.2)
$$
with $\ell +2m=n$. So $\Lambda_+$ is the sum of eigenspaces associated with
$\lambda_j=a_j$,
($j=1, \cdots, \ell$ )  $\lambda_j=c_{j-\ell}\pm id_{j-\ell}$, 
($j=\ell+1,\cdots,\ell+m$), 
and $\Lambda_-$ is the sum of eigenspaces associated with 
the corresponding $-\lambda_j$, and
$\Lambda_+ \oplus \Lambda_- = T_{\rho_0}{\bf R}^{2n}$. 
In these symplectic coordinates $\Lambda_+= \{ \xi=0 \},
\Lambda_-=\{
x=0\}$, and $F_{\rho_0}$ has block diagonal form,
the diagonal terms $\bigl( \lambda_1, \cdots, \lambda_\ell\bigr)$, the 
$2 \times 2$ matrices
$\normalbaselineskip=18pt \pmatrix {c_j&\pm d_j\cr \mp d_j & c_j\cr }$
($j=\ell+1,
\cdots, \ell+m$), 
the diagonal terms $\bigl( -\lambda_1, \cdots, -\lambda_\ell\bigr)$, and the 
$2 \times 2$ matrices
$\normalbaselineskip=18pt \pmatrix {-c_j&\mp d_j\cr \pm d_j & -c_j\cr }$
($j=\ell+1,
\cdots, \ell+m$) respectively, which is the so-called Cartan decomposition.
Note that $\Lambda_+$ and $\Lambda_-$ are dual spaces for the 
symplectic form on ${\bf R}^{2n}$.
To simplify notations, we shall sometimes 
introduce complex symplectic coordinates
$$\eqalign{
&z_{\ell+2j}={1\over\sqrt 2}(x_{\ell+2j}+ix_{\ell+2j-1}), \quad 
\zeta_{\ell+2j}={1\over\sqrt 2}(\xi_{\ell+2j}-i\xi_{\ell+2j-1}), \cr
&z_{\ell+2j-1}={1\over\sqrt 2}(x_{\ell+2j}-ix_{\ell+2j-1}),\quad 
\zeta_{\ell+2j-1}={1\over\sqrt 2}(\xi_{\ell+2j}+i\xi_{\ell+2j-1}), 
\quad j=1, \cdots, m\cr
}\leqno(1.3)$$
(the variables $x_j$ and $\xi_j$ being as in (1.2). ) Further we denote
$x_j$ for $z_j$, $\xi_j$ for the dual coordinate $\zeta_j$, and
eventually label
the collection of these symplectic coordinates, so that :
$$H_{p_2}= 
\Sum_{j=1}^n
\lambda_j (x_j {\partial
\over \partial x_j} - \xi_j {\partial
\over \partial \xi_j})\leqno(1.4)$$
or 
$$p''_{\rho_0}(t, s)=\Sum_{j=1}^n\lambda_j\bigl(t_{x_j}s_{\xi_j}+
t_{\xi_j}s_{x_j}\bigr)$$
Of course, we shall keep in mind that
the complexification here is only formal, since no analyticity
is assumed; this is no more than the usual
identification consisting for instance in taking complex coordinates
which diagonalize a rotation in the plane.

Now we turn to the non-linear case and 
recall the stable-unstable manifold
theorem. This theorem has a long history~: see e.g. 
[Ha] in the differentiable case, [Ch] or [Ne] for a proof
based on Sternberg's linearization theorem,
[AbMa], [AbRo] and references therein for more general
statements. Note that these results are generally
stated without symplectic structure, but most of them
easily extend to this setting. See however 
[Sj2,App] in the 
analytic category, and Theorem 2.2 below
for the almost holomorphic case.

\medskip
\noindent {\bf Theorem 1.1}: With notations above,
in a neighborhood of $\rho_0$, there are $H_p$-invariant lagrangian
manifolds ${\cal J}_\pm$ passing through $\rho_0$, such that
$T_{\rho_0}({\cal J}_\pm) = \Lambda_\pm$. Within ${\cal J}_+$ (resp. ${\cal
J}_-$), $\rho_0$ is repulsive (resp. attractive) for $H_p$, and $p|_{{\cal
J}_\pm}= 0$. We can also find real symplectic coordinates,
denoted again by $(x, \xi)$, such that their differential at $\rho_0$
verifies
$d(x,\xi)(\rho_0) = \id $, and 
${\cal J}_+= \{ \xi=0 \}, \ {\cal J}_-=\{ x=0\}$. 
In these coordinates
$$
p(x, \xi) = \langle A(x, \xi)x, \xi \rangle \leqno(1.5)
$$
where $A(x, \xi)$ is a real, $n \times n$ matrix with $C^\infty$
coefficients,
$A_0 = dA(\rho_0) =  \diag (\lambda_1, \cdots, \lambda_n)$
with the convention that if $\lambda_j$ is complex, diag$(\lambda_j, \overline
\lambda_j)$ denotes $\normalbaselineskip=18pt 
\pmatrix {c_j& -d_j\cr d_j & c_j\cr }$.
\medskip
It follows that
$$
H_p = A_1(x, \xi) x\cdot {\partial \over \partial x} - A_2(x, \xi) \xi
\cdot
{\partial
\over \partial \xi} 
\leqno(1.6)$$
with $A_j(x, \xi)=A_0+{\cal O}\bigl(x,\xi\bigr) $, 
$A_0=\diag (\lambda_1, \cdots, \lambda_n)$,
$A_1(x, \xi)=A(x, \xi)+{^t \hskip -1pt \partial}_\xi A(x, \xi)\cdot \xi$,
$A_2(x, \xi)={^t \hskip -1pt A}(x, \xi)+\partial_x A(x, \xi)\cdot x$,
and Spec $A(x, \xi)=$ Spec ${^t \hskip -1pt A}(x, \xi) \subset {\bf R}^+$.
Possibly after relabelling the coordinates, we may assume $ 0<\re
\lambda_1 \leq
\cdots \leq  \re \lambda_n$.

Now we describe the flow of $H_p$, using Proposition 1.1.
Let $\|\cdot\|$ denote the usual euclidean norm on ${\bf R}^n$. We put
$$B_0=\int_0^\infty e^{-s{^t \hskip -1pt A}_0}e^{-sA_0}ds$$
which is a positive definite
symmetric matrix, with the property
$^tA_0B_0+B_0A_0=I$. In the present case where $A_0$ is diagonalizable,
$$B_0=\diag \bigl(\lambda_1, \cdots, \lambda_l, {1\over2}c_{\ell+1},
{1\over2}c_{\ell+1}, \cdots, {1\over2}c_{\ell+m}, {1\over2}c_{\ell+m}\bigr)$$
If $\|x\|_0^2=\langle B_0x, x\rangle$ is the corresponding norm, then
$$A_0x\cdot\partial_x\|x\|_0^2=\|x\|^2, \quad 
A_0\xi\cdot\partial_\xi\|\xi\|_0^2=\|\xi\|^2 \leqno(1.7)$$
It follows from this and (1.6) that if $\|x\|_0^2+\|\xi\|_0^2\leq\delta^2$,
for some $\delta>0$ small enough, then
$${d\over dt}\|x\|_0^2=H_p\|x\|_0^2\geq C\|x\|^2, \quad 
-H_p\|\xi\|_0^2\geq C\|\xi\|^2 , \quad C>0$$
For $\delta >0$, we define
the outgoing region
$$\Omega_\delta^{\out }= \{ (x, \xi): \|\xi\|_0 < 2\|x\|_0, \
\|x\|_0^2+\|\xi\|_0^2 < \delta^2\} \leqno(1.9)$$ and let $\partial
\Omega_\delta^{\out }$ denote its boundary.
Let $t \mapsto (x(t), \xi(t)) = \exp tH_p \bigl(x(0), \xi(0)\bigr)$ be
an integral curve of $H_p$ with
$\rho=\bigl(x(0), \xi(0)\bigr) \in
\Omega_\delta^{\out }$. We have
$$\dot x(t) = A_1(x(t), \xi(t)) x(t), \ \dot \xi(t)=-A_2(x(t), \xi(t))\xi(t)
\leqno (1.10)$$
So when $\rho\in\Omega_\delta^{\out }$, $\|x(t)\|_0$ is increasing and 
$\|\xi(t)\|_0$ decreasing as long as $\bigl(x(t), \xi(t)\bigr)\in
\Omega_\delta^{\out }$, and moreover 
there is $C>0$ such that
for
$\delta>0$ sufficiently small and all $t\in{\bf R}$ :
$$\leqalignno{
&e^{-\re \lambda_+(t) t }e^{-C\delta |t|}  \|\xi(0)\|_0 \leq \|\xi(t)\|_0 \leq
e^{-\re \lambda_-(t)t}e^{C\delta |t| }
\|\xi(0)\|_0 &(1.11)\cr
&e^{\re \lambda_-(t)t}e^{-C\delta |t|} \|x(0)\|_0 \leq
\|x(t)\|_0 \leq e^{\re \lambda_+(t)t}e^{C\delta |t|}
\|x(0)\|_0 &(1.12)\cr 
}$$
with the convention $\lambda_+(t)=\lambda_n$ and $\lambda_-(t)=\lambda_1$
for $t>0$, $\lambda_+(t)=\lambda_1$ and $\lambda_-(t)=\lambda_n$
for $t<0$.
It follows that for any $\delta_0>0$, there is $\delta_1>0$ 
(say $\delta_1=\delta_0/2$, ) such that
if
$\rho \in \Omega_{\delta_1}^{\out }$, then exp$(-tH_p)(\rho)
\in
\Omega_{\delta_0}^{\out }$, $t \geq 0$, until the path meets
$\partial \Omega_{\delta_0}^{\out } \cap \{ \|\xi\|_0 = 2\|x\|_0 \}$.
For each $\rho \in \Omega_{\delta_1}^{\out }$, we define the hitting
time
$$T_-^{\out }(\rho) = \inf \{ t>0 : \|\xi(-t)\|_0 \geq 2 \|x(-t)\|_0 \},
\leqno(1.13)$$
i.e. the time for the path exp$(-tH_p)(\rho)$ to reach the cone
$\|\xi\|_0 = 2 \|x\|_0$. Since exp$(-tH_p)( \rho)$ is a $C^\infty$ function
of $\rho$ and $t$, it follows from the implicit function theorem
that $T_-^{\out }(\rho)$ is a $C^\infty$ function of $\rho \in
\Omega_{\delta_1}^{\out } \setminus {\cal J}_+$. For $\rho=(x, 0) \in
{\cal
J}_+$, we set
$T_-^{\out }(\rho)=+\infty$, and we leave it undefined for $\rho=0$.
Similarly, for $\rho \in \Omega_{\delta_1}^{\out }$ we define
$$T_+^{\out }(\rho) = \inf \{ t>0 : \|x(t)\|_0^2+\|\xi(t)\|_0^2 \geq
\delta_0^2
\},
\leqno(1.14)$$
to be the time for the path exp $(tH_p)(\rho)$ to leave the ball
$\|x\|_0^2+\|\xi\|_0^2 < \delta_0^2$. Again, $T_+^{\out }(\rho)$ is a
$C^\infty$ function of $\rho \in \Omega_{\delta_1}^{\out }$. Moreover,
there is $\tau>0$ such that for all $\rho \in \Omega_{\delta_1}^{\out
}$,
$\exp (tH_p)(\rho) \notin \Omega_{\delta_0}^{\out }$ for $T_-^{\out }(\rho)
\leq t \leq T_-^{\out }(\rho) + \tau$.
Since we are interested in local properties of the flow near
$\rho_0$, we can modify, without loss
of generality, $p(x, \xi)$ outside a small neighborhood of $\rho_0$
such that the path exp $(tH_p)(\rho)$, $\rho \in \Omega_{\delta_1}^{\out }$,
will
never enter again $\Omega_{\delta_0}^{\out }$ after time  $T_+^{\out
}(\rho)$,
i.e. we may assume $\tau=+\infty$. From now on, we change notation
$\delta_0$ and $\delta_1$ to $\delta$ for simplicity, keeping in mind
that $\delta$ is a sufficiently small, but fixed positive number.

We define in a similar way the incoming region
$$\Omega_\delta^{\i }= \{ (x, \xi): \|x\|_0 < 2\|\xi\|_0, \
\|x\|_0^2+\|\xi\|_0^2 < \delta^2\} \leqno(1.15)$$
and the hitting times $T_\pm^{\i }(\rho)$. More precisely, 
$$\leqalignno{
&T_-^{\i }(\rho) = \inf \{ t>0 : \|x(-t)\|_0^2+\|\xi(-t)\|_0^2 \geq
\delta^2 \} &(1.16) \cr
&T_+^{\i }(\rho) = \inf \{ t>0 : \|x(t)\|_0\geq 2 \|\xi(t)\|_0 \}
&(1.17)\cr
}$$
As above, we may assume that the flow starting from any point 
$\rho\in{\bf R}^{2n}$ crosses at most once the region $\Omega_\delta
=\Omega_\delta^{\i }
\cup \Omega_\delta^{\out }$. Then estimates (1.11) and (1.12)
hold for all $(x, \xi)\in
\Omega_\delta$, and all $t\in{\bf R}$ provided 
$\bigl(x(t), \xi(t)\bigr)\in
\Omega_\delta$.
\medskip
Now let $I$ denote the ideal of $C^\infty({\bf R}^{2n})$ 
consisting in all smooth functions vanishing at $\rho_0$. We want to
solve the homological equation $H_pf=g$ in
$I^\infty$. This is of course essentially wellknown~: see e.g. 
[GuSt,p.175] for analogous results. So
let $\chi^{\out }+\chi^{\i }=1$ be
a smooth partition of unity in the unit sphere ${\bf S}^{2n-1}$ such
that
supp $\chi^{\out} \subset \{ \|\xi\|_0 < 2\|x\|_0 \}$,
supp $\chi^{\i } \subset \{ \|x\|_0 < 2\|\xi\|_0 \}$. We extend
$\chi^{\out }, \ \chi^{\i }$ as homogeneous functions of 
degree 0 on $T^*{\bf R}^n\setminus\rho_0$.

\medskip
\noindent {\bf Proposition 1.2}~: Let $\rho_0$ be an hyperbolic fixed
point for
$p$ as above, and $g
\in I^\infty$. Let 
$$
f^{\out }(\rho) = \int_{-\infty}^0 \bigl( \chi^{\out }g \bigr)\circ \exp
(tH_p)(\rho) dt, \ f^{\i }(\rho) = -\int_0^\infty  \bigl( \chi^{\i }g
\bigr)\circ \exp (tH_p)(\rho) dt
$$
Then  
$f = f^{\out }+f^{\i }\in I^\infty$ solves $H_pf=g$.
\medskip
\noindent {\it Proof}: We treat the case of $f^{\out }$, this of $f^{\i }$
is similar. Let $\delta_0>0$ small enough, and 
$\Omega_{\delta_1}^{\out / \i }$ be as above. Without loss of generality,
we may assume supp $g\subset \Omega_{\delta_0}= \Omega_{\delta_0}^{\out } 
\cup  \Omega_{\delta_0}^{\i }$, so $\supp (\chi^{\out } g)\subset  
\Omega_{\delta_0}^{\out }$. Then it is easy to see that 
$$(\supp f^{\out})\cap \Omega_{\delta_1}\subset  \Omega_{\delta_1}^{\out }$$
so we will assume $\rho\in \Omega_{\delta_1}^{\out }$,
and as above write $\delta$ for $\delta_0$ or $\delta_1$. If 
$\rho\in \Omega_{\delta}^{\out }\setminus {\cal J}_+$, we have
$f^{\out }(\rho) = \int_{-T_-^{\out }(\rho)}^0 
\bigl( \chi^{\out }g \bigr)\circ \exp
(tH_p)(\rho) dt$, since $\exp
(tH_p)(\rho) \notin \supp \chi^{\out }$ for $t<-T_-^{\out }(\rho)$.
Furthermore, 
$$H_pf^{\out }(\rho)=\int_{-\infty}^0 
{d\over dt}\bigl(\bigl( \chi^{\out }g \bigr)\circ \exp
(tH_p)(\rho)\bigr) dt=(\chi^{\out }g)(\rho)$$ 
When $\rho\in {\cal J}_+$, $\exp(tH_p)(\rho)\to 0$ when $t\to -\infty$
and the integral makes sense because of (1.12) and the fact that 
$g(\rho)={\cal O}(\rho)$, as $\rho\to 0$. Again  
$H_pf^{\out }(\rho)=\chi^{\out }g(\rho)$. We are left to show that 
$f^{\out } \in I^\infty$. Because of (1.12) and $\|\xi(t)\|_0\leq
2\|x(t)\|_0$ in supp $\chi^{\out }$, $f^{\out }$ is continuous and
vanishes at $\rho=0$. To show that $f^{\out } \in C^1$, we write,
following [IaSj]~:
$$d\bigl( (\chi^{\out }g)\circ \exp(tH_p)(\rho)\bigr)
=\bigl(d (\chi^{\out }g)( \exp(tH_p)(\rho)\bigr)\circ d\exp(tH_p)(\rho)
\leqno(1.19)$$
so we need to examine the evolution of $d\kappa_t(\rho)=d\exp(tH_p)(\rho)$
along the integral curve $\kappa_t$ of $H_p$ starting at $\rho$.
Differentiating $\partial_t\kappa_t(\rho)=H_p\bigl(\kappa_t(\rho)\bigr)$
we find 
$$\partial_td\kappa_t(\rho)={\partial H_p\over \partial \rho}
(\kappa_t(\rho))\circ \bigl(d\kappa_t(\rho)\bigr), \quad d\kappa_0(\rho)=\id
\leqno(1.20)$$
with ${\partial H_p\over \partial \rho}(\rho)=2F_{\rho_0}+{\cal O}(\rho)$,
and Gronwall lemma applied to (1.20), as in (1.11) and (1.12) gives for 
$\kappa_t(\rho)\in\Omega_\delta^{\out }$, all $t\leq 0$:
$$\leqalignno{
&e^{-(\re \lambda_1 -C\delta)t }  \leq \|d{\xi}_t(\rho)\| \leq
e^{-(\re \lambda_n + C\delta)t }
&(1.21)\cr
&e^{(\re \lambda_n+C\delta)t}  \leq
\|dx_t(\rho)\| \leq e^{(\re \lambda_1-C\delta)t}
&(1.22)\cr 
}$$
so $d\kappa_t(\rho)={\cal O}\bigl(e^{-(\re \lambda_n + C\delta)t }\bigr)$.

On the other hand, $g$ being flat at 0, 
$d (\chi^{\out }g)( \exp(tH_p)(\rho)= {\cal O}(\|x_t(\rho)\|^N)$,
any $N$, so taking $N$ large enough, we see that 
$d\bigl( (\chi^{\out }g)\circ \exp(tH_p)(\rho)\bigr)$ is integrable,
so $f^{\out }\in C^1$, and vanishes at 0. To continue, we take partial
derivative of (1.20) with respect to $\rho_j$, $j=1, \cdots, 2n$ and write
$$
\partial_t{\partial \over \partial \rho_j}
d\kappa_t(\rho)-{\partial H_p\over \partial \rho}
(\kappa_t(\rho))\circ \bigl({\partial \over \partial \rho_j}
(d\kappa_t(\rho)\bigr)=F_j(t, \rho)
\leqno(1.23)$$
with 
$$F_j(t, \rho)=\Sum_{k=1}^{2n}{\partial^2H_p \over \partial \rho_k
\partial \rho}(\kappa_t(\rho)){\partial \over \partial \rho_j}
\kappa_{t, k}(\rho)\circ d\kappa_t(\rho)\leqno(1.24)$$
Using the group property, we write (1.20) as 
$$\partial_td\kappa_{t-\widetilde t}\bigl(\kappa_{\widetilde t}(\rho)\bigr)
\circ d\kappa_{\widetilde t}(\rho)
={\partial H_p\over \partial \rho}
\bigl(\kappa_{t-\widetilde t}\bigl(\kappa_{\widetilde t}(\rho)\bigr)
\bigr)\circ 
d\kappa_{t-\widetilde t}\bigl(\kappa_{\widetilde t}(\rho)\bigr)
\circ d\kappa_{\widetilde t}(\rho), \quad d\kappa_0(\rho)=\id
$$
Since $\kappa_{\widetilde t}$ is a canonical map, 
$d\kappa_{\widetilde t}$ is invertible, so
$$\partial_td\kappa_{t-\widetilde t}\bigl(\kappa_{\widetilde t}(\rho)\bigr)
={\partial H_p\over \partial \rho}
\bigl(\kappa_{t-\widetilde t}\bigl(\kappa_{\widetilde t}(\rho)\bigr)
\bigr)\circ 
d\kappa_{t-\widetilde t}\bigl(\kappa_{\widetilde t}(\rho)\bigr), 
\quad d\kappa_0(\rho)=\id
$$
So we recognize $d\kappa_{t-\widetilde t}\bigl(\kappa_{\widetilde t}(\rho)
\bigr)$, $d\kappa_{\widetilde t-\widetilde t}(\rho)=\id $ as the fundamental
matrix of our $2n \times 2n$ system of ordinary differential equations,
and Duhamel's principle gives, since 
${\partial \over \partial \rho_j}
d\kappa_t(\rho)|_{t=0}=0$ :
$${\partial \over \partial \rho_j}
d\kappa_t(\rho)=\int_0^t d\kappa_{t-\widetilde t}
\bigl(\kappa_{\widetilde t}(\rho)\bigr)\circ
F_j(\widetilde t, \rho)d\widetilde t\leqno(1.25)$$
Fom (1.21) and (1.22) we find the estimate 
$F_j(\widetilde t, \rho)={\cal O}\bigl(e^{-2(\re \lambda_n + C\delta)
\widetilde t }\bigr)$, and by integration
$${\partial \over \partial \rho_j}
d\kappa_t(\rho)=
{\cal O}\bigl(e^{-2(\re \lambda_n + C\delta)t }\bigr)\leqno(1.26)$$
On the other hand, differentiating (1.19) with respect to $\rho_j$ we get :
$$\eqalign{
{\partial \over \partial \rho_j}
d\bigl( (\chi^{\out }g)&\circ\kappa_t(\rho)\bigr)=
d (\chi^{\out }g)\bigl( \kappa_t(\rho)\bigr)\circ 
{\partial \over \partial \rho_j}d\kappa_t(\rho)+\cr
&+\Sum_{k=1}^{2n}{\partial \over \partial \rho_k}
d\bigl( \chi^{\out }g\bigr)\bigl( \kappa_t(\rho)\bigr)
{\partial \over \partial \rho_j}\kappa_{t,k}(\rho)
\circ d\kappa_t(\rho)\cr}\leqno(1.27)
$$
Using (1.26), and again (1.21), (1.22), the estimates 
$$d (\chi^{\out }g)\circ\bigl(\kappa_t(\rho)\bigr),\quad
{\partial \over \partial \rho_k}
d\bigl( \chi^{\out }g\bigr)\bigl( \kappa_t(\rho)\bigr)
= {\cal O}\bigl( \|x_t(\rho)\|^N\bigr)$$
ensure once more the integrability of
${\partial \over \partial \rho_j}
d\bigl( (\chi^{\out }g)\circ\kappa_t(\rho)\bigr)$,
so $f^{\out }\in C^2$ and we can see that its second derivatives vanish
at 0. The argument carries over easily by induction,
so the Proposition is proved. $\clubsuit$ 
\medskip
Note that we used here for convenience $C^\infty$ coordinates adapted to
${\cal J}_\pm$, but the proof is essentially independent of coordinates
(see a variant of this in Proposition 4.3 below.)
\medskip
Now we are ready for proving Theorem 0.1, by combining the Birkhoff
normal form (see e.g. Appendix for a simple proof)
and a deformation argument. When $p$ has a non degenerate
critical point with non-resonant frequencies, we know that there is
a smooth canonical transform $\kappa$ between neighborhoods of 0,
leaving fixed the origin, such that $p\circ\kappa(x, \xi)=q_0(\iota)
+r(x, \xi)$, where $\iota=(\iota_1,\cdots,\iota_n)$ are the action
variables as in (0.3), and $r\in I^\infty$ depends also on 
the corresponding dual (angle) variables. The hamiltonian $q_0(\iota)$
satisfies the same hypotheses as $p$, and is constructed from the 
formal Taylor series by a Borel sum of the type
$q_0(\iota)=\Sum_{k=1}^\infty \widetilde q_k(\iota)\chi\bigl(
{\iota\over {\e }_k}\bigr)$, $\chi\in C^\infty_0({\bf R}^n)$
equal to 1 near 0, ${\e }_k\to 0$ fast enough as $k\to \infty$,
and $\widetilde q_k(\iota)$ is homogeneous of degree $k$.
The canonical transformation is of the form
$\kappa=\exp H_{\widetilde f}$
for some smooth $\widetilde f$. We shall try to construct a family
$\kappa_s$ of canonical transformations, $0\leq s\leq 1$,
tangent to identity at infinite order, such that $\kappa_0=I$
and $\kappa_1$ solves $p\circ\kappa\circ \kappa_1=q_0$. The deformation
(or homotopy) method consists in finding a $C^\infty$ field $X_s$
along which some property is conserved,
in that case the property for a smooth family of hamiltonians, interpolating
between $p$ and $q_0$, of being integrable.
It reduces here essentially to
solving a homological equation as in Proposition 1.2
(see [ArVaGo] for an introduction, and also [GuSc,p.168], [MeSj],
[BaLlWa], [IaSj] \dots, for other applications 
more directly relevant to our problem.) 
So let $q_s=q_0+sr$, $0\leq s\leq 1$, and look for $\kappa_s$ such that
$$q_s\circ\kappa_s=q_0\leqno(1.30)$$
Then $\kappa_s|_{s=1}$ will solve our problem. The deformation field 
$$X_s(\rho)=\Sum_{j=1}^{2n}v_{s,j}(\rho){\partial \over \partial\rho_j}\in 
I^\infty(T{\bf R}^{2n})$$
is such that 
$$\partial_s\kappa_s=v_s\circ\kappa_s \leqno(1.31)$$
Differentiating (1.30) gives 
$r\circ\kappa_s+{\partial q_s\over \partial \rho}
(\kappa_s)\circ \partial_s\kappa_s=0$, or
$$
r\circ\kappa_s+\langle v_s\bigl(\kappa_s(\rho)\bigr), q_s\bigl(
\kappa_s(\rho)\bigr)\rangle=0$$
Furthermore, we require $X_s$ to be hamiltonian, i.e.
$v_s=H_{f_s}$, $f_s\in I^\infty$, so we get
$$\langle H_{f_s}, q_s\rangle=-\langle H_{q_s}, f_s\rangle=-r \leqno(1.32)$$
all quantities being evaluated at $\kappa_s(\rho)$. 
We want to apply Proposition 1.2 to $p=q_s$, $g=r$, so we move to the
new symplectic coordinates (adapted to the outgoing/incoming manifolds)
by composing with smooth canonical  
transformations $\Phi_s$, i.e. replace $H_{q_s}$ by $(\Phi_s)^*H_{q_s}$, 
$f_s$ by $(\Phi_s)^*f_s$, etc\dots, so omitting for brevity these
coordinates transformations when no confusion might occur, Proposition 1.2
gives $f_s\in I^\infty$
solving (1.32). So we are led to show,
that given $H_{f_s}\in I^\infty$, (1.31) has a solution of the 
form $\kappa_s=I+\kappa'_s$, $\kappa'_s\in I^\infty$.
Existence for $0\leq s\leq 1$
follows e.g. from Gronwall lemma, truncating $q_s$ outside a neighborhood
of 0,
and the condition $\kappa_0=I$
gives 
$$\|\kappa_s(\rho)\|\leq C\|\rho\|, \quad C>0 \leqno(1.33)$$ 
for $\|\rho\|<\delta$.
We want to show 
$\kappa'_s(\rho)={\cal O}(\rho^\infty)$. Recall from the proof of 
Proposition 1.2 that, by the group property, $d\kappa_s(\rho)$
is the fundamental solution for the system 
$\partial_sY(\rho, s)={\partial H_{f_s}\over \partial\rho}
\bigl(\kappa_s(\rho)\bigr)Y(\rho, s)$. 
Since $d\kappa'_s(\rho)$ solves 
$$\partial_s
d\kappa'_s(\rho)-{\partial H_{f_s}\over \partial \rho}
(\kappa_s(\rho))\circ \bigl(
d\kappa'_s(\rho)\bigr)={\partial H_{f_s}\over \partial \rho}
(\kappa_s(\rho)), \quad d\kappa'_s(0)=0
\leqno(1.34)$$
Duhamel's principle gives 
$$
d\kappa'_s(\rho)=\int_0^s d\kappa_{s-\widetilde s}
\bigl(\kappa_{\widetilde s}(\rho)\bigr)\circ
 {\partial H_{f_s}\over \partial \rho}
(\kappa_{\widetilde s}(\rho)) d\widetilde s
$$
Since ${\partial H_{f_s}\over \partial \rho}
(\kappa_{\widetilde s}(\rho))={\cal O}(\|(\kappa_{\widetilde s}(\rho))\|^N)$,
(1.33) gives ${\partial H_{f_s}\over \partial \rho}
(\kappa_{\widetilde s}(\rho))={\cal O}(\|\rho\|^N)$,
and 

$d\kappa_{s-\widetilde s}
\bigl(\kappa_{\widetilde s}(\rho)\bigr)={\cal O}(1)$,
so chosing $N$ large enough, we get
$d\kappa'_s(\rho)={\cal O}(\|\rho\|^2)$.
Integrating this relation, we get again $\kappa'_s(\rho)={\cal O}(\|\rho\|)$.
Taking partial derivative of (1.34) with respect to $\rho_j$ as in the
proof of Proposition 1.2 yields also ${\partial\over\partial\rho_j}
d\kappa'_s(\rho)={\cal O}(\|\rho\|)$, and a straightforward induction
argument shows $\kappa'_s\in I^\infty$, uniformly for $s$
on compact sets. Taking $s=1$ 
and undoing the transformation $\Phi_s|_{s=1}$
give eventually the result. $\clubsuit$

\medskip
We pause for a while, presenting our result in some different way.
It is sometimes convenient to perform the 
Birkhoff transform in action-angle
coordinates (see [Ga,p.473] for the elliptic case.~) We restrict
for simplicity to the usual case of a (real-) hyperbolic fixed point, where  
$$p(x, \xi)= \xi^2-\Sum_{j=1}^n \lambda_j^2 x_j^2
+{\cal O}(\|x\|^3)$$ 
The corresponding Williamson coordinates are then given by the
linear symplectic transformation 
$\kappa_1(x, \xi)= (y, \eta)$, $\sqrt 2\lambda_jy_j=\lambda_jx_j+\xi_j$,
$\sqrt 2\lambda_j\eta_j=-\lambda_jx_j+\xi_j$. We define hyperbolic
action-angle coordinates $(\iota, \varphi)$ by the formulas
$\lambda_jx_j=\sqrt {2\lambda_j\iota_j}\cosh \varphi_j$, 
$\xi_j=\sqrt {2\lambda_j\iota_j}\sinh \varphi_j$, and set
$\kappa_0(\iota, \varphi)=(x, \xi)$. Let $\kappa$
be the canonical transform of theorem 0.1, and define 
$\widetilde \kappa= \kappa_0^{-1}\circ\kappa_1^{-1}\circ\kappa\circ
\kappa_1\circ\kappa_0$. 
Then, with $\kappa(y, \eta)=(y', \eta')=(y, \eta)+{\cal O}(|y, \eta|^2)$,
we have $\widetilde \kappa(\iota, \varphi)=(\iota', \varphi')$,
$2\lambda_j\iota^\prime_j-2\lambda_jy^\prime_j\eta_j=-{\xi^\prime_j}^2+\lambda_j^2{x^\prime_j}^2$,
where $\kappa_1(x', \xi')=(y', \eta')$. 
There exists a smooth
generating function 
$S(\iota', \phi)$  such that 
$\iota=\partial_\varphi S(\iota', \varphi)$, $\varphi'=\partial_{\iota'}
S(\iota', \varphi)$, and of the form $S(\iota', \varphi)=\langle\iota',
\varphi\rangle+
\Phi(\iota',\varphi)$. Here $\partial_{\iota'}\Phi(\iota',\varphi)=
{\cal O}(\iota')$, $\partial_\varphi\Phi(\iota',\varphi)=
{\cal O}(\iota'^2)$, uniformly for $\varphi$ in compact sets, and
$\iota'$ small enough. Moreover, $p=q(\iota')$. 

\medskip
\noindent {\bf b) Integrability near a closed trajectory of hyperbolic 
type}.
\smallskip
In this section we consider an hamiltonian flow with a non trivial
center manifold. More precisely, let $p=p_E$ be a smooth, real 
(family of) hamiltonian(s)
on ${\bf R}^{2n}$ ($E$ is one of the $2n$ variables,~) 
and $K$ the set of trapped trajectories near 
energy 0~:
$$K=\{\rho\in p^{-1}_E(0), \ E\in J=[-\e _0, \e_0], \
\exp (tH_{p_E})(\rho)\not\to \infty,\
\hbox{as} \ t\to \pm\infty \}$$
Let $K_{\e }=K\cap p^{-1}_E(0)$, $E=\e $ small, and assume for simplicity we 
are in the situation where 
$K_0=\gamma_0$ 
is a closed 
trajectory of hyperbolic type. This is the case when $p_E$ is a function of
$2(n-1)$
phase variables $(x', \xi')\in T{\bf R}^{n-1}$, periodic
with respect to $\theta\in{\bf S}^1$~; parameter $E$ then stands for the dual
variable. 

Then in a neighborhood of $K$, there is 
a smooth, symplectic, closed submanifold $\Sigma$ of dimension 2,
containing $K_0$ and such that $H_{p_E}$ is tangent to $\Sigma$
everywhere. We call $\Sigma$ the center manifold of $\gamma_0$,
and it is nothing but the one-parameter family of closed trajectories
$\gamma_{\e }\subset p^{-1}_E(0)$, $E=\e $ small. Hyperbolicity means that 
$p_E$ vanishes of second order on $\Sigma$, and for all $\rho\in\Sigma$,
the fundamental matrix $F_{\rho}$ as in (1.1) is of rank $2n-2$,
and has no purely imaginary eigenvalues. In the case at hand, we will
assume that these eigenvalues are rationally independent.
For $\rho\in\Sigma$, 
let as above $\Lambda_{\pm}(\rho)\subset T_\rho({\bf R}^{2n})$
be the $(n-1)$-dimensional isotropic subspaces whose complexifications
are the sum of all complex eigenspaces corresponding to eigenvalues
with positive/negative real parts. We have the splitting
$(T_\rho\Sigma)^\perp=\Lambda_+(\rho)\oplus\Lambda_-(\rho)$,
where $(\bullet)^\perp$ stands for `` symplectic orthogonal''.
The restriction $\sigma_\Sigma$ of $\sigma$ to $T\Sigma^\perp$
is clearly invariant under $H_{p_E}$. 
Again, we recall the center-stable-unstable manifold theorem extending
Theorem 1.1~:
\medskip
\noindent {\bf Theorem 1.6}: With notations above,
in a neighborhood of $\Sigma$, there are (unique)
$H_{p_E}$-invariant, smooth involutive
manifolds ${\cal J}_\pm$ passing through $\Sigma$, 
such that
for all $\rho\in\Sigma$,
$T_{\rho}({\cal J}_\pm) = \Lambda_\pm(\rho)$. Within ${\cal J}_+$ 
(resp. ${\cal J}_-$), $\Sigma$ is repulsive (resp. attractive) for $H_{p_E}$, 
and $p_E|_{{\cal
J}_\pm}= 0$ (recall that $E$ is one of the variables.~)
We can also find real symplectic coordinates,
denoted again by $(x, \xi)=\bigl((x', x''), (\xi', \xi'')\bigr)$,
such that their differential verifies
$d(x,\xi)|_\Sigma = \id $,
$\Sigma$ is given by $(x', \xi')=0$,
${\cal J}_+= \{ \xi'=0 \}$ and ${\cal J}_-=\{ x'=0\}$. 
In these coordinates
$$
p_E(x, \xi) = \langle A(x, \xi)x', \xi' \rangle \leqno(1.40)
$$
where $A(x, \xi)$ is a real, $(n-1)\times (n-1)$ matrix with $C^\infty$
coefficients, and eigenvalues
$\lambda_1(x'',\xi''), \cdots, \lambda_{n-1}(x'',\xi'')$.
\medskip
Of course, ${\cal J}_+$ depend on $E$, and also on $\theta$ 
that we have omitted 
in the notations.  We may now also forget the variable $E$. 
Theorem 1.6  is proved e.g. as in Theorem 2.2 below. 

Our constructions extend readily to this situation. We still define the 
outgoing/incoming region, for instance 
$$\Omega_f^{\out }= \{ (x, \xi): \|\xi'\|_0 < 2\|x'\|_0, \
\|x'\|_0^2+\|\xi'\|_0^2 < f(x'',\xi'')\} \leqno(1.41)$$ 
where $f$ is a smooth, positive function with sufficiently small support and 
small derivatives. 

Now let $I_\Sigma$ denote the ideal of $C^\infty({\bf R}^{2n})$ 
consisting in all smooth functions in $\{\|x'\|_0^2+\|\xi'\|_0^2 < f(x'',\xi'')
\}$ vanishing at $\Sigma$. We choose
as above 
a smooth partition of unity $\chi^{\out }+\chi^{\i }=1$ in the unit sphere 
${\bf S}^{2n-3}$ such
that
supp $\chi^{\out} \subset \{ \|\xi'\|_0 < 2\|x'\|_0 \}$,
supp $\chi^{\i } \subset \{ \|x'\|_0 < 2\|\xi'\|_0 \}$, and extend
$\chi^{\out }, \ \chi^{\i }$ as homogeneous functions of 
degree 0 on $T^*{\bf R}^n\setminus\Sigma$. Then for
$p$ as above, and $g \in I_\Sigma^\infty$, if  
$$
f^{\out }(\rho) = \int_{-\infty}^0 \bigl( \chi^{\out }g \bigr)\circ \exp
(tH_p)(\rho) dt, \ f^{\i }(\rho) = -\int_0^\infty  \bigl( \chi^{\i }g
\bigr)\circ \exp (tH_p)(\rho) dt
$$
then  
$f = f^{\out }+f^{\i }\in I_\Sigma^\infty$ solves $H_pf=g$.

Let $\rho_0\in\Sigma$ be such that the non resonance condition holds
on the eigenvalues 
$\lambda_1(\rho_0)$, \dots, $\lambda_{n-1}(\rho_0)$, and
apply the Birkhoff normal form to $p$. 
Then there exist a smooth canonical transform $\kappa$ for the symplectic
2-form $\sigma_\Sigma$, 
and a smooth hamiltonian $q_0(\iota')$, where $\iota'=(\iota_1, \cdots,
\iota_{n-1})$ are action variables as in (1.3) built from the 
$(x',\xi')$-coordinates,
such that
$$p\circ \kappa(x, \xi)=q_0(\iota')+r(x, \xi), \ r\in I_\Sigma^\infty,
\ (x, \xi)\in \neigh \ (\rho_0, {\bf R}^{2n})
\leqno (1.42)$$
Next we pass to the deformation procedure, composing with a new canonical
transformation, preserving $\sigma_\Sigma$, to remove the remainder
$r$. So we get, with a new $\kappa$~: 
$$p\circ \kappa(x, \xi)=q_0(\iota'), \ (x, \xi)\in \neigh 
(\rho_0, {\bf R}^{2n}) \leqno(1.43)$$
To formulate a semi-global result we assume that the fundamental 
matrix of $p$ (for the 2-form $\sigma_\Sigma$) is constant on $\Sigma$, 
with non resonant frequencies as above. The constructions above depending
smoothly on $\rho_0\in\Sigma$, we have found a smooth fibre bundle
over $\Sigma$, foliated by action-angle coordinates in $T\Sigma^\perp$
adapted to $p$. The question of triviality for this bundle is left
open. See [CuB], [Vu] for other
(semi-)global aspects of integrability.

\medskip
\noindent {\bf 2. The Lewis-Sternberg normal form for the Poincar\'e map}.
\medskip
In this section we prove Theorem 0.2. First we recall the following version
of a theorem of Lewis-Sternberg ([St,Thm1,Corollary1.1], [Fr,ThmV.1] and
[IaSj] for a detailed proof. ) For simplicity we content to a particular 
case relevant to
our problem.  
So assume $A$ is a real $2n\times 2n$ symplectic matrix and has eigenvalues 
$\lambda_1, \cdots, \lambda_n$, $1/\lambda_1,
\cdots, 1/\lambda_n$, $\overline\lambda_1,
\cdots, \overline\lambda_n$, $1/\overline\lambda_1,
\cdots, 1/\overline\lambda_n$, 
none of them is of modulus 1. Then there is a natural choice of the 
logarithm 
$B=\log A$, and $B$ is antisymmetric for the canonical 2-form on 
$T^*{\bf R}^n$. Let $\mu_j=\log \lambda_j$, in such a way that $\overline 
\lambda_j$ corresponds to $\overline \mu_j$, and $p_0(\rho)=b(\rho)={1\over 2}
\sigma(\rho, B\rho)$. 
Assume that 
for $k_j\in{\bf Z}$ 
$$\Sum k_j\mu_j\in 2i\pi{\bf Z} \Longrightarrow \Sum k_j\mu_j= 0 \leqno(2.1)$$
We have the following :
\medskip
\noindent {\bf Theorem 2.1}: Let 
$\Phi : \neigh (0, {\bf R}^{2n}) \to \neigh (0, {\bf R}^{2n})$
be a smooth canonical transformation, leaving fixed $\rho_0=0$,
and $A=d\Phi(\rho_0)$ as above. Then there  
$p\in C^\infty$ defined near $\rho_0$, uniquely determined modulo $I^\infty$,
(for a given choice of $p_0$ ) such that $p(\rho)=p_0(\rho)+{\cal O}(\rho^3)$
and 
$$\Phi(\rho)=\exp H_p(\rho)+{\cal O}(\rho^\infty) \leqno (2.2)$$
\medskip
Let us sketch from [IaSj] the main ideas of the proof. As above, is relies
on a deformation argument. Given $p_s$ a smooth real function depending
smoothly on the real parameter $s$, $p_s(\rho)=b(\rho)+{\cal O}(\rho^3)$,
we consider the corresponding
canonical transformation 
$$\Phi_{t,s}=\exp tH_{p_s}$$
Since $p_s$ vanishes to second order at $\rho_0$, the germ of $\Phi_{t,s}$
at $\rho_0$ is well defined for all real $t$.
Arguing as in Proposition 1.2, the first variation with respect to $s$ is
integrated between 0 and $t$, which gives~: 
$$\partial_s \Phi_{t,s}
=(\Phi_{t,s})_* H_{q_{t,s}}\leqno(2.3)$$ 
where 
$$q_{t,s}=\int_0^t\partial_s p_s
\circ \Phi_{\widetilde t,s} d\widetilde t \leqno(2.4)$$
In this formula, we take $t=1$ (deleting the 
corresponding subscript) and consider a problem where 
$\partial_s p_s$ will be the unknown. More precisely, starting from 
$\Phi_0=\exp H_b$, we want to find $p$ such that 
$\Phi(\rho)= \exp H_p(\rho)+{\cal O}(\rho^\infty)$, by interpolating
with a family 
$$\Phi_s(\rho)= \exp H_{p_s}(\rho)+{\cal O}(\rho^\infty), \ 0\leq s\leq 1$$
If we define $q_s(\rho)={\cal O}(\rho^3)$
by $\partial_s \Phi_s =(\Phi_s)_* H_{q_s}$ as in (2.3), then we solve 
(2.4) by successive approximations, as
$$q_s=\int_0^1\partial_s p_s^{(N)}
\circ \exp tH_{p_s^{(N)}}(\rho) dt+{\cal O}(\rho^{N+1}) \leqno(2.5)$$  
with $p_s^{(0)}=b(\rho)$. This can be done precisely because of
hypothesis (2.1), and it can be shown that the sequence
$p_s^{(N)}$ is asymptotic, for $\rho$ near 0, to a $C^\infty$ function
$p_s(\rho)$. Taking $s=1$ gives (2.2) with $p=p_1$. 
Uniqueness follows again from a deformation argument. $\clubsuit $
\medskip
Assume further as in Theorem 0.2 that the $\mu_j$'s are rationally
independent, i.e. 
$$\Sum k_j\mu_j=0 \Longrightarrow k_j= 0 \leqno(2.6)$$
Using Birkhoff transformations shows that there is a smooth $\kappa$,
$\kappa(\rho_0)=\rho_0$ 
such that $p\circ\kappa(\rho)=q(\rho)+{\cal O}(\rho^\infty)$ and $q=q(\iota)$ 
depends on the action
variable only. So we have 
$$\kappa^{-1}\circ \exp H_p \circ \kappa = \exp H_q \ \hbox{mod}\ I^\infty
\leqno(2.7)$$
and again by Theorem 2.1 
$$\kappa^{-1}\circ \Phi \circ \kappa = \exp H_q \ \hbox{mod}\ I^\infty
\leqno(2.8)$$
Recall the following result, which is the symplectic version
of the Lewis-Sternberg theorem, and 
whose proof is very similar 
to our Theorem 0.1~:
\medskip 
\noindent {\bf Theorem 2.2} [BaLLWa]: Let 
$f, g : \neigh (0, {\bf R}^{2n}) \to \neigh (0, {\bf R}^{2n})$
be smooth canonical transformations, leaving fixed $\rho_0=0$,
and assume they are tangent to infinite order at $\rho_0$. 
Let $A=df(\rho_0)$ have its spectrum outside the unit
circle as above. Then there is a smooth canonical transform $\chi$ 
leaving fixed $\rho_0$, $d\chi(\rho_0)=\id $, such that 
$\chi^{-1}\circ f \circ \chi =g$. 
\medskip
Now it is clear that Theorem 0.2 immediately follows from Theorem 2.2 
and (2.7), (2.8) applied to 
$f=\kappa^{-1}\circ \Phi \circ \kappa$, and $g=\exp H_q$.
$\clubsuit$

\medskip
\noindent {\bf 3. Parameter dependent case}.
\medskip
We extend some of the previous results, taking advantage of the 
fact observed in [IaSj], that the Birkhoff normal form can be 
carried out nearby critical points with non resonant frequencies.
\medskip
\noindent {\bf a) The Birkhoff normal form}.  
\smallskip
Let $p^{(s)}\in C^\infty$ as in the Appendix, depend 
smoothly on $s\in \neigh (0, {\bf R}^k)$, 
$p^{(s)}(\rho_0)=0$, and have a non-degenerate critical point 
of hyperbolic type at $\rho_0$. (In some applications, the critical
point depends on $s$, but choosing suitable linear
symplectic coordinates
and changing $p^{(s)}$ by a constant we are in this situation.~)  
Possibly after
performing another linear symplectic transformation, we may assume
that its quadratic part
is of the form 
$$p^{(s)}_2(x, \xi)=\Sum_{j=1}^n\mu_j^{(s)}x_j\xi_j \leqno (3.1)$$
with coordinates independent of $s$. 
(For simplicity we take real frequencies. )
For $s=0$, we suppose
the $\mu_j=\mu_j^{(0)}$ rationally independent. 
Then Proposition A.1 below shows there is a smooth family
of canonical transforms, $\kappa^{(s)}$, $\kappa^{(s)}(\rho_0)=\rho_0$,
such that 
$$p^{(s)}\circ\kappa^{(s)}(\rho)=q^{(s)}(\iota)+r^{(s)}(\rho), \
r^{(s)}(\rho)={\cal O}(\rho^\infty)+
{\cal O}(s^\infty\rho^3)
\leqno (3.2)$$
with the principal part of $q^{(s)}$ as in (3.1). 
Looking at the deformation procedure, we see that we can apply
the stable/unstable manifold theorem to  
$Q_\sigma(\rho)=q^{(s)}(\iota)+\sigma r^{(s)}(\rho)$,
$0\leq\sigma\leq 1$, and if we decompose 
$r^{(s)}=u^{(s)}+v^{(s)}$, $v^{(s)}={\cal O}(\rho^\infty)$, 
$u^{(s)}={\cal O}(s^\infty\rho^3)$, we are able to solve 
$H_{Q_\sigma}f_\sigma=v^{(s)}$,
for $f_\sigma \in I^\infty$.
Then the vector field $X_\sigma=H_{f_\sigma}$
generates a 1-parameter family of canonical transformations
$\kappa_\sigma$ as in (1.31), and for $\sigma=1$ we get 
$$p^{(s)}\circ\kappa^{(s)}\circ\kappa_1(\rho)=q^{(s)}(\iota)+
{\cal O}(s^\infty\rho^3)\leqno(3.3)$$
which is the normal form for $p^{(s)}$. 
\medskip
\noindent {\bf b) The Lewis-Sternberg normal form}.
\smallskip
As in [IaSj] we extend Theorem 0.2 to the parameter dependent case,
thinking for instance of the Poincar\'e map that depends smoothly
on the energy. 
For simplicity we just vary one parameter $s\in\neigh (0, {\bf R})$.
Let 
$\Phi^s : \neigh (0, {\bf R}^{2n}) \to \neigh (0, {\bf R}^{2n})$,
$s\in\neigh (0, {\bf R})$,
be a smooth family of smooth canonical transformations, 
leaving fixed $\rho_0=0$,
and $A^s=d\Phi^s(\rho_0)$. We assume that $\Phi=\Phi^0$
fulfills the assumptions of Theorem 2.1, so that 
$$\Phi(\rho)=\exp H_p(\rho)+{\cal O}(\rho^\infty)$$
where $p=p^0$ is unique modulo $I^\infty$ and 
the choice of its quadratic part. 
As above, we want to extend $p$ to a smooth real-valued family
$p^s$, with~:
$$\Phi^s(\rho)=\exp H_{p^s}(\rho)+\rho^2{\cal O}\bigl((\rho, s)^\infty\bigr)
\leqno(3.5)$$
Define $q^s={\cal O}(\rho^2)$ as in (2.4) (with $t=1$) so that 
$(\Phi^s)^*\partial_s \Phi^s= H_{q^s}$. At the level of
Taylor expansions, we replace (2.4) by its approximation
and the problem is to find $q^s$ such that~:
$$q^s=\int_0^1\partial_s p^s
\circ \exp tH_{p^s}(\rho) dt+\rho^2{\cal O}\bigl((\rho, s)^\infty\bigr)
\leqno(3.6)$$  
We try to achieve this condition at any order in $s$.
At zeroth order, i.e. for $s=0$, we get a unique solution 
$\partial_s p^s|_{s=0}={\cal O}(\rho^2)$, mod $I^\infty$. 
If we differentiate $k$ times we get
$$\int_0^1(\partial_s^{k+1} p^s)
\circ \exp tH_{p^s}(\rho) dt=\partial_s^kq^s(\rho)+
F_k(p^s, \cdots, \partial_s^kp^s, \rho)+
\rho^2{\cal O}\bigl((\rho, s)^\infty\bigr)$$  
and if $p^0, \cdots, \partial_s^kp^s|_{s=0}={\cal O}(\rho^2)$
have been determined, we get $\partial_s^{k+1} p^s|_{s=0}={\cal O}(\rho^2)$
from this equation. It is then clear that (3.6) has a solution which
is unique modulo $\rho^2{\cal O}\bigl((\rho, s)^\infty\bigr)$. Let
$\widetilde \Phi^s=\exp H_{p^s}$. Then 
$$(\Phi^s)^*\partial_s \Phi^s=(\widetilde\Phi^s)^*\partial_s \widetilde\Phi^s
+\rho^2{\cal O}\bigl((\rho, s)^\infty\bigr), \widetilde\Phi^0=\Phi^0$$
and it easily follows that (3.5) holds. 

Assume now that for $s=0$
the $\mu_j$'s are rationally
independent. Using the parameter dependent Birkhoff transformations 
as in Proposition A.1, we see
that for $s\in\neigh (0, {\bf R})$ small enough, there is a smooth 
family of hamiltonians $q^s$, and 
canonical transformations $\kappa^s$,
$\kappa^s(\rho_0)=\rho_0$,
such that $p\circ\kappa^s=q^s+{\cal O}(\rho^\infty)+
{\cal O}(\rho^3 s^\infty)$ and $q^s=q^s(\iota)$ 
depend on the action
variable only. So we have 
$$(\kappa^s)^{-1}\circ \exp H_{p^s} \circ \kappa^s = 
\exp H_{q^s} +{\cal O}(\rho^\infty)+{\cal O}(\rho^3 s^\infty)
\leqno(3.7)$$
and by (3.5) : 
$$(\kappa^s)^{-1}\circ \Phi^s \circ \kappa^s = 
\exp H_{q^s} + \rho^2{\cal O}\bigl((\rho, s)^\infty\bigr)
\leqno(3.8)$$
There suffices then to apply a parameter dependent version of 
Theorem 2.2 (which follows easily from a careful inspection of the proof
in [BaLLWa],~)
to see that (3.7) and (3.8) imply the following~:
\medskip
\noindent {\bf Proposition 3.1}: Let $\Phi^s$, $s\in \neigh (0, {\bf R})$, 
be a smooth family of smooth canonical transformations,
$\Phi^s: \neigh (0, {\bf R}^{2n})\to \neigh (0, {\bf R}^{2n})$, 
$\Phi^s(0)=0$, such that for $s=0$, $A^0=d\Phi^0(0)$
is non degenerate, its eigenvalues $\lambda_j$, $j=1\cdots,n$,
lie outside the unit circle,
and $\mu_j=\log \lambda_j$ verify (2.1). Assume further that 
the $\lambda_j$'s are rationally independent.
Then there are 
a smooth family of smooth canonical maps $\kappa^s$, 
$s\in \neigh (0, {\bf R})$,
$\kappa^s(0)=(0)$, $d\kappa^s(0)=\id $, and   
a smooth family of smooth functions $q^s(\iota)$ 
depending on the action variables 
$\iota$ alone, such that 
$$(\kappa^s)^{-1}\circ \Phi^s \circ \kappa^s = 
\exp H_{q^s}+\rho^2{\cal O}(s^\infty)$$

\medskip
\noindent {\bf 4. The complex case}.
\medskip
We present here a rather rough discussion in the almost holomorphic case,
i.e. for hamiltonians whose $\overline \partial$ vanishes of infinite order
at $\rho_0$, somewhat in the spirit of [Sj2,3] and [MeSj].
First we recall some properties concerning symplectic structures in 
$T{\bf C}^n$~; then we state the center stable/unstable manifold theorem
for almost holomorphic hamiltonians~; at last we 
prove Theorem 0.4, and conclude with some elementary properties on
monodromy.
\medskip
\noindent {\bf a) Complex symplectic geometries}.
\medskip
The variables in the complex phase-space
$T^*{\bf C}^n$ will still be denoted by $(x, \xi)$. As in the real case,
we start with some geometric preparations.

First we recall some elementary facts about complex vector fields. If 
$$v(\rho)=\Sum_{j=1}^{2n}v_j(\rho)\partial_{\rho_j}+
v'_j(\rho)\overline \partial_{\rho_j}
=\Sum_{j=1}^{2n}\bigl( a_j(\rho)\partial_{x_j}+b_j(\rho)\partial_{\xi_j}
+a'_j(\rho)\overline\partial_{x_j}+b'_j(\rho)\overline \partial_{\xi_j}
\bigr) \in T(T^*{\bf C}^n)$$ 
is a vector field on $T^*{\bf C}^n$, 
we set $\widehat v=2\re v= v+\overline v$, or
$$\widehat v(\rho)=\Sum_{j=1}^{2n}
\bigl(v_j(\rho)+\overline {v'_j(\rho)}\bigr)\partial_{\rho_j}
+\bigl(\overline {v_j(\rho)}+v'_j(\rho)\bigr) \  \overline \partial_{\rho_j}
\leqno(4.1)$$
Identifying ${\bf C}^n\times{\bf C}^n$ with ${\bf R}^{2n}\times{\bf R}^{2n}$, 
$\widehat v$ is simply the vector 
$$\eqalign{
(v_1+&\overline{v'_1}, \cdots, v_{2n}+\overline{v'_{2n}})
=(a_1+\overline{a'_1}, \cdots, a_n+\overline{a'_n},
b_1+\overline{b'_1}, \cdots, b_n+\overline{b'_n})=\cr
&=\bigr(\re (a_1+a'_1), \im (a_1-a'_1), \cdots,\re (b_n+b'_n), \im (b_n-b'_n)
\bigr)\cr
}$$
expressed in the basis 
$B=\bigl(\partial_{\re x_1}, \partial_{\im x_1},
\cdots, \partial_{\re \xi_n}, \partial_{\im \xi_n}\bigr)$.  
In general the identification between ${\bf C}^n$ (or ${\bf C}^{2n}$)
and the underlying real vector space will be expressed as
$\Theta (a_1,\cdots,a_n)=(\re a_1, \im a_1, \cdots, \re a_n, \im a_n)$.

Let us denote by $I$ the ideal of
$C^\infty$ functions in ${\bf C}^n$ (or $T^*{\bf C}^n$ as will be clear
from the context, ) that vanish at $\rho_0$. 
We assume throughout that $v'_j\in I$, or even $v'_j\in I^\infty$. 
In that case, we write $v\in T^{(1,0)}(T^*{\bf C}^n)\oplus
T^{(0,1)}_\infty(T^*{\bf C}^n)$.
Then $\widehat v$ is the (unique) real vector field which gives the same 
result 
as $v$, at the point $\rho_0$, when
applied to a differentiable function $u$,  provided
$\overline \partial u \in I$.
For real $t$, the flow of $\widehat v$
will be denoted by 
$$\widehat \Phi_t(\rho)=\bigl(\widehat x_t(\rho), \widehat \xi_t(\rho)\bigr)=
\exp (t\widehat v)(\rho)$$
In the case where $v'_j=0$ (i.e. $v\in T^{(1,0)}(T^*{\bf C}^n)$, )
this is the solution of the system of ODE's
$${d\over dt}\widehat {(x_j)}_t(\rho)=a_j\bigl(\widehat \Phi_t(\rho)\bigr), 
\quad {d\over dt}\widehat {(\xi_j)}_t(\rho)=b_j
\bigl(\widehat \Phi_t(\rho)\bigr), \quad \widehat \Phi_0(\rho)=\rho$$ 
So it has the property, that if $v\in T^{(1,0)}(T^*{\bf C}^n)$ 
has holomorphic coefficients, then
$\widehat \Phi_t(\rho)$ is the restriction to the real $t$-axis
of the holomorphic flow 
$$ \Phi_t(\rho)=\bigl( x_t(\rho),  \xi_t(\rho)\bigr)=
\exp (t v)(\rho) \leqno (4.2)$$
We recall also  that  ${\bf C}^{2n}$
is endowed with the complex canonical 2-form $\sigma_{\bf C}$,
which makes it a symplectic space, and 2 real symplectic 2-forms 
$\re \sigma_{\bf C}$ and
$\im \sigma_{\bf C}$. For convenience, we remove 
subscript ${\bf C}$ from the notations.
If $p$ is a smooth complex function on ${\bf C}^{2n}$, the hamiltonian vector 
field of $p$ is defined as
$$H_p={\partial p\over \partial \xi}{\partial\over \partial x}+
{\partial p\over \partial \overline \xi}{\partial\over \partial \overline x}
-{\partial p\over \partial x}{\partial\over \partial \xi}-
{\partial p\over \partial \overline x}{\partial\over \partial \overline \xi}
$$
(note we have used a different convention as in [MeSj], [Sj1],
where $H_p$ does not contain the antiholomorphic derivatives.)
If we define the real hamiltonian 
vector field $H^{\re \sigma}$ by 
$(\re \sigma) (H^{\re \sigma}_f, t)=\langle df, t\rangle$, then we have
$H^{\re \sigma}_{\re p}=\widehat H_p$. More precisely,
in the basis $B$~:
$$\widehat H_p=\bigl( \re {\partial p \over \partial\re \xi},
-\re {\partial p \over \partial\im \xi},
-\re {\partial p \over \partial\re x},
\re {\partial p \over \partial\im x}\bigr)\leqno(4.3)$$
We denote by ${\partial \widehat H_p \over \partial\rho}$
the Jacobian (in the real sense) expressed in this basis.

The Poisson bracket associated with $\re \sigma_{\bf C}$ is denoted
by $\{\bullet, \bullet \}_R$ and coincides with $\{\re \bullet, \re \bullet\}$
for the real symplectic structure on ${\bf C}^{2n}$ read through 
$\Theta$. 

Let $p$ be a smooth function such that $\overline \partial p \in I^\infty$.
For real $t$, the hamiltonian flow of $\widehat H_p$
will be denoted as above by 
$$\widehat\Phi_t(\rho)=\bigl(\Phi_{t,x}(\rho), \Phi_{t,\xi}(\rho)\bigr)=
\exp (t\widehat H_p)(\rho) \leqno (4.4)$$
Let $\widetilde X_\rho=\Theta\bigl(\overline \partial_x\Phi_{t,x}, 
\overline \partial_x\Phi_{t,\xi}\bigr)$, and 
$\widetilde Y_\rho=\Theta\bigl(\overline \partial_\xi\Phi_{t,x}, 
\overline \partial_\xi\Phi_{t,\xi}\bigr)$ considered as vector fields
on $T^*({\bf C}^n)$. In the same way, we write
$ X_\rho=\Theta\bigl( \partial_x\Phi_{t,x}, 
\partial_x\Phi_{t,\xi}\bigr)$, and 
$Y_\rho=\Theta\bigl( \partial_\xi\Phi_{t,x}, 
 \partial_\xi\Phi_{t,\xi}\bigr)$, where $\partial$ denotes the holomorphic
derivative.
We state first some technical Lemma, which follows from a staightforward
computation, and the fact that $p$ verifies approximately 
the Cauchy-Riemann equations~:
\medskip

\noindent {\bf Lemma 3.1}: With $p$ as above
$$\leqalignno{
&\partial_t\widetilde X_\rho
={\partial \widehat H_p \over \partial\rho}(\widehat{\Phi_t})
\widetilde X_\rho+
{\cal O}\bigl(\widehat\Phi_t(\rho)^\infty\bigr)\bigl(  
\widetilde X_\rho,  X_\rho
\bigr) &(4.5)\cr
&\partial_t\widetilde Y_\rho
={\partial \widehat H_p \over \partial\rho}(\widehat\Phi_t)
\widetilde Y_\rho+
{\cal O}\bigl(\widehat\Phi_t(\rho)^\infty\bigr)\bigl(  
\widetilde Y_\rho, Y_\rho
\bigr) &(4.6)\cr
}$$
\medskip
\noindent {\bf b) The stable-unstable-center manifold
theorem in the complex domain}.
\medskip
Our first step is to extend the stable/unstable manifold theorem
in the case of almost holomorphic hamiltonians. To be complete we will 
actually prove 
a little bit more than required.
We will follow closely the 
nice geometric argument
of [Sj2] in the analytic category, implemented for higher derivatives
by idea we borrowed also from [HeSj1].  

So let $p$ such that $\overline \partial p \in I^\infty$, have a 
non degenerate critical point at $\rho_0$, $p(\rho_0)=0$.
Let $F_{\rho_0}(p)$ as in (1.1)
denote the fundamental matrix (in the holomorphic sense,~) and  
assume as before that $F_{\rho_0}(p)$ has $2n$ distinct, no purely 
imaginary eigenvalues. Let again $\Lambda_\pm \subset T_{\rho_0}{\bf C}^{2n}$ 
be the sum of all eigenspaces
corresponding to eigenvalues with positive (resp. negative) real
parts. We have~:
\medskip

\noindent {\bf Theorem 4.2} : With the notations above,
in a neighborhood of $\rho_0$, there are 
$\widehat H_p$-invariant, R-lagrangian
manifolds ${\cal J}_\pm$ (i.e. lagrangian for $\re \sigma_{\bf C}$~,)
passing through $\rho_0$, such that
$T_{\rho_0}({\cal J}_\pm) = \Lambda_\pm$. Within ${\cal J}_+$ (resp. ${\cal
J}_-$), $\rho_0$ is repulsive (resp. attractive) for $\widehat H_p$, 
and $\re p|_{{\cal
J}_\pm}= 0$. We can also find $\re \sigma_{\bf C}$-symplectic coordinates,
denoted again by $(x, \xi)=\kappa (y, \eta)$, 
$\overline \partial \kappa \in I^\infty$, 
such that their differential at $\rho_0$
verifies
$d\kappa(\rho_0) = \id $, 
$\kappa^*(\sigma_{\bf C})=\sigma_{\bf C}$ mod $I^\infty$
and ${\cal J}_+= \{ \xi=0 \}, \ {\cal J}_-=\{ x=0\}$. 
In these coordinates
$$
\re p(x, \xi) = \re \langle A(x, \xi)x, \xi \rangle 
\leqno(4.8)$$
where $A(x, \xi)$ is smooth, has constant term equal to $A_0$,
and $\overline\partial A(x, \xi)\in 
I^\infty$. Moreover, 
$$p(x, \xi) = \langle A(x, \xi)x, \xi \rangle \ \hbox{mod } \ I^\infty
\leqno(4.9)$$

(For simplicity,
we have written $\langle A(x, \xi)x, \xi \rangle$ instead of 
$\langle A'(x, \xi)(\re x,\im x), (\re \xi, \im \xi) \rangle$
where $A'(x, \xi)$ is a $2n \times 2n$ matrix~; actually the notation
$p(x, \xi)=\langle A(x, \xi)x, \xi \rangle$ makes sense
at the level of formal Taylor expansion at $\rho_0$.~) 

\medskip

\noindent {\it Outline of the Proof}: 
We proceed in several steps ; in the topological
step we start 
to define, as in Sect.1.a, 
the outgoing/incoming regions  
relative to $\widehat H_p$, and study the flow of
lagrangian manifolds, as $t\to\pm\infty$. This yields, via
a compactness argument, to
$C^0$ coordinates where the outgoing (resp. incoming) submanifold
${\cal J}_+$ (resp. ${\cal J}_-$ ) is given by $\xi=0$ (resp. $x=0$ )~;
then, we turn to differentiability and prove the ${\cal J}_+$
are $C^1$. Then we turn to higher derivatives and properties of 
almost analyticity.

We first choose coordinates where $F_{\rho_0}$ has block-diagonal form. 
Taking complex linear coordinates as in (1.3), we can make it
diagonal. Then the hamiltonian vector field takes the form
$$H_p=A_0x\cdot{\partial \over \partial x}-A_0\xi\cdot
{\partial \over \partial \xi}+{\cal O}\bigl( \|x,\xi\|^2 \bigr) 
\bigl({\partial \over \partial x}, {\partial \over \partial \xi}\bigr) 
\mod \ T^{(0,1)}_\infty(T^*{\bf C}^n) \leqno(4.10)$$
where we recall $A_0=\diag (\lambda_1, \cdots, \lambda_n)$.
For real $t$, let $\widehat\Phi_t(\rho)$
be the hamiltonian flow of $\widehat H_p$ as in (4.4).
As in Sect.1.a we can construct an hermitian norm $\|\cdot\|_0$
such that identity (1.7) holds if $\|\cdot\|$
and $\|\cdot\|_0$ stand now for the hermitian norms.
For $\delta >0$, we define
the outgoing region as
$$\Omega_\delta^{\out }= \{ (x, \xi): \|\xi\|_0 < 2\|x\|_0, \
\|x\|_0^2+\|\xi\|_0^2 < \delta^2\}$$ 
and let $\partial
\Omega_\delta^{\out }$ denote its boundary.
Estimates (4.10)
again show that there exists $C>0$ such that~:
$$\|x\|_0/(2C)\leq\widehat H_p \|x\|_0\leq  \|x\|_0/C, 
\quad \rho \in \Omega_\delta^{\out }
\leqno(4.12)$$
while
$$-\widehat H_p \|\xi\|_0\geq  \|\xi\|_0/C \
\hbox{on} \ \partial\Omega_\delta^{\out }\cap 
\{ (x, \xi): \|\xi\|_0 = 2\|x\|_0\} \leqno(4.13)$$ 
Let $t \mapsto \widehat \Phi_t \bigl(x(0), \xi(0)\bigr)$ be
an integral curve of $\widehat H_p$ with
$\rho=\bigl(x(0), \xi(0)\bigr) \in
\Omega_\delta^{\out }$. 
Along $\widehat \Phi_t $, we have $\partial_t=\widehat H_p$,
so using (4.12), Gronwall Lemma, after suitably
truncating $p$ outside a neighborhood of $\rho_0$, shows that 
$$e^{t/(2C)}\|x(0)\|_0\leq \|x(t)\|_0\leq e^{t/C}\|x(0)\|_0, \ 
\rho \in \Omega_\delta^{\out }, t\geq 0
\leqno(4.14)$$
which allows to define the hitting times $T_\pm^{\out }$
as in (1.13) and (1.14)
(although we have not yet found the outgoing manifold.)

It follows from (4.14) and (4.13) that $\Omega_\delta^{\out }$
is stable under $\widehat \Phi_t$, i.e. if 
$\rho \in\Omega_\delta^{\out }$, then 
$\exp(t\widehat H_p)(\rho)\in\Omega_\delta^{\out }$, for
$0\leq t < T_\pm^{\out }(\rho)$, while it never gets back in afterwards.

Similarly, we define the incoming region as in (1.15), and the 
corresponding hitting times as in (1.16), (1.17).

Now we try to find the outgoing/incoming manifolds for $\widehat H_p$,
and study the evolution of the complex manifold 
$\Lambda_t=\{\exp(t\widehat H_p)
(\rho): \rho\in \Omega_\delta^{\out } \}$, as $ t\to +\infty$. It is convenient
to introduce 
$$\Lambda^{\out }=\{ \bigl(t, \tau; \exp(t\widehat H_p)
(\rho)\bigr), \rho\in\Omega_\delta^{\out },
0\leq t < T_+^{\out }(\rho), \tau=\re p\bigl(\exp(t\widehat H_p)(\rho)\bigr)
\}\leqno(4.15)$$
By what we have just said, $\Lambda^{\out }$ is a connected submanifold
of codimension 1 in the symplectic space $T^*{\bf R}^{2n+1}$ endowed with
the 2-form
$d\tau\wedge dt+\re \sigma_{\bf C}$. The vector
field $\partial_t+\widehat H_p$ is tangent to $\Lambda^{\out }$, and $\tau$
is independent of $t$. 
The evolution of a tangent vector 
$\widehat X_\rho (t)=\bigl( 
\widehat X_x(t), \widehat X_\xi(t) \bigr)\in T(T^*{\bf R}^{2n})$
(the $\rho$-projection of the tangent space to $\Lambda^{\out }$, )
is given by the $4n \times 4n$ system
$$\partial_t \widehat X_\rho(t)={\partial\widehat H_p\over \partial\rho}
\bigr(\widehat\Phi_t(\rho)\bigr) \widehat X_\rho (t)
\leqno(4.16)$$
where $\partial_\rho$ denotes the gradient in the real sense.

It is easy to see that the leading term in the $4n\times 4n$ matrix
${\partial\widehat H_p\over \partial\rho}$ in the basis $B$
has a hyperbolic structure, each eigenvalue $\lambda_j$ 
occuring twice, as well as $-\lambda_j$, $\pm \overline \lambda_j$, 
so that the linear flow is expansive in the $(\re x, \im x)$-
directions, and contractive in the $(\re \xi, \im \xi)$-
directions.  

So (4.16) shows that 
if $\e _0>0$ and $\delta>0$ are sufficiently small, then the 
outgoing region  
$$\|\widehat X_\xi(t)\|\leq \e _0 \|\widehat X_x(t)\| \leqno(4.17)$$
is stable along $\widehat \Phi_t(\rho)$, $\rho\in\Omega_\delta^{\out }$, 
as $t$ increases, $t<T^{\out }_+(\rho)$.

Now let $\widetilde{\cal J}_+ =\{\rho \in \Omega_\delta^{\out }; \xi=0\}$, 
${\cal J}_+(t)=\widehat \Phi_t(\widetilde {\cal J}_+)\cap
\Omega^{\out }_\delta$,
and $\Lambda_+$ be its lift in $\Lambda^{\out }$. This is a 
submanifold of $T^*{\bf R}^{2n+1}$, lagrangian for 
$d\tau\wedge dt+\re \sigma_{\bf C}$, and its tangent space contains
$\partial_t+\widehat H_p$. 
Applying the 
theorem of constant rank to the projection 
$\pi : \Lambda_+ \to {\bf C}^n_x$, (4.17) 
shows that $\Lambda_+ $
(or ${\cal J}_+(t)$, forgetting about $\tau$ which is independent of $t$, 
and that we may take equal to 0, since $p(\rho_0)=0$,~)
is of the form $\xi=g_+(t, x)$ where $g_+\in C^\infty$ (see for instance [M]
for a simple proof. ) 
Moreover, $g_+(0, x)=0$.   
Since $\widehat \Phi_t(\rho)\in \Omega_\delta^{\out }$, we have
$\|g_+(t, x)\|_0\leq 2\|x\|_0$ for all $t\geq 0$.
By compactness, there is a sequence
$t_j\to +\infty$, such that $g_+(t_j, \cdot)\to g_+$ in 
$C^0\bigl(\{\|x\|< \hbox{Const.} \delta\}\bigr)$. 
We put
$${\cal J}_+=\{(x, g_+(x)) : x \in \neigh (0)\} \leqno(4.19)$$
(the outgoing tail, or outgoing manifold) 
and proceed to show that $g_+\in C^1$. 

Consider the evolution 
of a normal vector 
$$\widehat Z_\rho (t)=\bigl( 
\widehat Z_x(t), \widehat Z_\xi(t) \bigr)\in N({\cal J}_+(t))=
\bigl( T(T^*{\bf R}^{2n})|_{{\cal J}_+(t)}\bigr) / T({\cal J}_+(t))
\leqno (4.21)$$
(the $\rho$-projection of the normal space to $\Lambda_+$. )
It is given by :
$\partial_t \widehat Z_\rho(t)= M(\widehat \Phi_t(\rho)) \widehat Z_\rho (t)$
where the leading part of $M(x, \xi)$ is obtained from this of 
${\partial\widehat H_p\over \partial\rho}$ by permuting the eigenvalues
with positive and negative real parts. So in $\Lambda_+$, the region given 
by 
$$\|\widehat Z_\xi(t)\|\geq  \|\widehat Z_x(t)\|/\e _0 \leqno(4.22)$$
is stable under $\widehat \Phi_t$. 

Let now $\rho_t$ be another integral curve of $\widehat H_p$,
starting at $\rho\in\Omega_\delta^{\out }$, and
not in ${\cal J}_+(t)$
($\rho_t$ lies in $\Lambda^{\out }$, but we choose the initial condition
away from $\widetilde {\cal J}_+$.~)
Let $\Gamma_t$ be the orthogonal
projection of $\rho_t$ on ${\cal J}_+(t)$, ${\dot \Gamma}_t\in
N({\cal J}_+(t))$ 
the normal
vector. By (4.22), we see that if $\gamma_t$ 
denotes the lenght of the segment $[\rho_t,\Gamma_t]$, 
then ${d\over dt}\gamma_t
\leq - C \gamma_t$, $C>0$~; so the integral curves of $\widehat H_p$
approach ${\cal J}_+$ exponentially fast as $t$ increases,
and the estimate $\|g_+(t, x)-g_+(s, x)\|={\cal O}(e^{-s/C})$,
all $t\geq s\geq 0$, shows that $g_+(t, x)$ is Cauchy,
and  $T({\cal J}_+(t))$ has a limit
as $t\to +\infty$ (not only for a subsequence $t_j$. ) This limit
is the tangent space to ${\cal J}_+=\{\xi=g_+(x) : x\in 
\neigh (0)\}$, and it follows that ${\cal J}_+(t)$
tends exponentially fast to ${\cal J}_+$ in the $C^1$ topology.
It is easy to see that ${\cal J}_+$ is invariant under $\widehat
\Phi_t$, all $t$, and characterized as the set of
$\rho\in\Omega_\delta^{\out }$ such that 
$\widehat\Phi_t(\rho)\in\Omega_\delta^{\out }$, all $t\leq 0$. We have
$\widehat \Phi_t(\rho)\to \rho_0=0$ as $t\to -\infty$, 
$\rho\in {\cal J}_+$. Moreover, $\re p=\tau=0$ on ${\cal J}_+$.

We are left to show that ${\cal J}_+$ is a lagrangian submanifold for
$(T^*{\bf C}^n, \re \sigma_{\bf C})$. If $u_1, \ u_2$ are complex $C^1$
functions vanishing on ${\cal J}_+$, and $\rho\in{\cal J}_+$, 
then $\{u_1, u_2\}_R(\rho)=\{u_1\circ \Phi_t, u_2\circ \Phi_t\}_R
(\Phi_{-t}(\rho))$. Since integral curves of $\exp t\widehat H_p$
approach ${\cal J}_+$ exponentially fast, we see that 
$du_j\circ \Phi_t
\bigl(\Phi_{-t}(\rho)\bigr)$ tends to 0 as $t\to +\infty$, hence
$\{u_1, u_2\}_R =0$, 
and we have proved that ${\cal J}_+$ is involutive.
Because $T_{\rho_t}{\cal J}_+(t)$ is transversal to
$\widetilde {\cal J}_-=\{\rho\in \Omega^{\out }_\delta, x=0\}$ 
(and their intersection is 0)
we have also proved, letting $t\to +\infty$,
that ${\cal J}_+$ is lagrangian for $\re\sigma_{\bf C}$.
Furthermore,
$T_{\rho_0}({\cal J}_+)=\Lambda_+$. 
Similarly, we introduce
$$\Lambda^{\i }=\{ \bigl(t, \tau; \exp(-t\widehat H_p)
(\rho)\bigr), \rho\in\Omega_\delta^{\i },
0\leq t < T_-^{\i }(\rho), \tau=\re p\bigl(\exp(t\widehat H_p)(\rho)\bigr)
\}\leqno(4.24)$$
Taking the flow of $\widetilde {\cal J}_-$ 
through 
$\widehat\Phi(t)$ for negative $t$, we set    
${\cal J}_-(t)=\widehat \Phi_t(\widetilde {\cal J}_-)\cap
\Omega^{\i }_\delta$, and look for the evolution of a tangent vector
to ${\cal J}_-(t)$ along an integral curve $\rho_t$ of $\widehat H_p$,
starting at $\rho\in \Omega^{\i }_\delta$, and not in ${\cal J}_-(t)$. 
Letting $t\to -\infty$, we can see that
${\cal J}_-(t)$ tends exponentially fast
to ${\cal J}_-=\{(g_-(\xi), \xi) : \xi\in \neigh (0)\}$, for some $C^1$
function $g_-(\xi)$. Then ${\cal J}_-$  is again 
lagrangian with respect to $\re \sigma_{\bf C}$, and we call it
the incoming tail, or incoming manifold. Again we have $\re p=\tau=0$
on ${\cal J}_-$.

It is clear that the invariant manifolds ${\cal J}_\pm$ are characterized 
as the set of $\rho\in\Omega_\delta$ such that 
$\widehat \Phi_{\mp t}(\rho)\in\Omega_\delta$, for all $\pm t\geq 0$.

The higher derivatives cannot apparently be handled with the same method,
but by the uniqueness property of the outgoing/incoming manifolds,
we can conclude as in [AbRo,App.C] with a fixed point argument, the 
limits being necessarily ${\cal J}_\pm$. An alternative way is to follow
the proof of [HeSj1,Prop.2.3].
Namely, it follows easily from the previous arguments
that ${\cal J}_+(t)$ (say) can be parametrized by a phase function
$\varphi_t(x, \eta)$, such that the graph of
$\exp(t\widehat H_p)$, $t\geq 0$,  is given by
$$C_t=\{\bigl( \partial_\eta\varphi_t, \eta, x, \partial_x\varphi_t\bigr),
\ (x, \eta)\in \Omega_\delta^{\out } \}$$
Furthermore, $\varphi_t$ verifies the eikonal equation 
$${\partial \varphi_t\over \partial t}+\re p\bigl(x, 
{\partial \varphi_t\over \partial x}\bigr)=0, \ \varphi|_{t=0}=
\langle x, \eta \rangle
$$
By the previous estimates, we know then that $\varphi_t$ tends 
exponentially fast as $t \to +\infty$, to some $\varphi_+(x, \eta)$ 
in $C^2(\Omega_\delta^{\out })$. Then $\varphi_+(x, \eta)$ verifies
again the corresponding stationary eikonal equation, and parametrizes
${\cal J}_+$.  
Using the transport
equations verified by ${\partial \varphi_t\over \partial x}$,
we can show as in [HeSj1] that this convergence holds actually in 
$C^\infty(\Omega_\delta^{\out })$. We proceed similarly in 
$\Omega_\delta^{\i }$.

Once we have found the smooth, involutive invariant manifolds ${\cal J}_\pm$,
we choose adapted coordinates of the 
form $(x', \xi')=\bigl(x-g_-(\xi), \xi-g_+(x)\bigr)$. 
By construction, these are smooth
symplectic coordinates for $\re \sigma_{\bf C}$,
where the outgoing (resp. incoming) manifold takes the form
$\xi'=0$ (resp. $x'=0$.) From now on, we work in these coordinates,
which we denote again by $(x, \xi)$, deleting 
the prime. The same argument as in Sect.1
then shows that (1.11) and (1.12) hold for $\rho\in\Omega_\delta$,
$t \in {\bf R}$, where 
$(x(t), \xi(t))$ stands for $\widehat\Phi_t(\rho)$, and $\|\cdot\|_0$ for the
hermitean norm.

We pass now to the almost analyticity property.
Using coordinates adapted to ${\cal J}_\pm$,
this can be done again by combining Lemma 4.1 with 
the method above, showing that the generating functions verify 
$\overline \partial \varphi_\pm \in I^\infty$. (Alternatively, this can be
done by the fixed point argument of [AbRo,App.C]. )
The Theorem easily follows, since also (4.9) can be 
recovered  from (4.8), using that $p$ verifies the Cauchy-Riemann
equations modulo $I^\infty$. 
$\clubsuit$

\medskip
\noindent {\bf c) Proof of Theorem 0.4}.
\medskip
We proceed exactly as in the
real case. Let again $\chi^{\out }+\chi^{\i }=1$ be 
a smooth partition of unity in $T^*{\bf R}^{2n}\setminus\rho_0$
with
supp $\chi^{\out} \subset \{ \|\xi\|_0 < 2\|x\|_0 \}$,
supp $\chi^{\i } \subset \{ \|x\|_0 < 2\|\xi\|_0 \}$. We start with the~:
\medskip
\noindent {\bf Proposition 4.3}: Let $p$ be as above, and $g\in I^\infty$. Let 
$$
f^{\out }(\rho) = \int_{-\infty}^0 \bigl( \chi^{\out }g \bigr)\circ \exp
(t\widehat H_p)(\rho) dt, \ f^{\i }(\rho) = -\int_0^\infty  \bigl( \chi^{\i }g
\bigr)\circ \exp (t\widehat H_p)(\rho) dt
$$
Then  
$f = f^{\out }+f^{\i }\in I^\infty$ solves $\widehat H_pf=g$.
\medskip
(Note that we may avoid using the $C^\infty$ coordinates where
${\cal J}_\pm$ are given by $\xi=0$ and $x=0$, as we did in 
Proposition 1.2, since the proof is essentially independent of the choice of
coordinates. All what really matters is the existence of the 
differentiable manifolds ${\cal J}_\pm$.~)  

Using again Birkhoff series (in ${\bf C}^{2n}$), we know that there is
a smooth canonical transform for the complex symplectic structure 
$(T^*{\bf C}^n, \sigma_{\bf C})$,
$\kappa(\rho_0)=\rho_0$, and
such that 
$$p\circ\kappa(x, \xi)=q_0(\iota)+r(x, \xi)\leqno(4.27)$$ 
where $\iota=(\iota_1,\cdots,\iota_n)$ are the action
variables as in (0.3), and $r\in I^\infty$. 
The hamiltonian $q_0(\iota)$
satisfies the same hypotheses as $p$, and is constructed from the 
formal Taylor series by a Borel sum of the type
$q_0(\iota)=\Sum_{k=1}^\infty \widetilde q_k(\iota)\chi\bigl(
{\iota\over {\e }_k}\bigr)$, $\chi\in C^\infty_0({\bf C}^n)$
equal to 1 near 0, of the form $\chi(z_1, \cdots, z_n)=\chi_0(z_1)
\otimes \cdots \otimes \chi_0(z_n)$, $\chi_0$ rotation invariant.
Of course, $\overline \partial_\iota q_0(\iota)={\cal O}(\iota^\infty)$.
Using again Borel sums, the canonical transformation is of the form
$\kappa=\exp H_{\widetilde f}$
for some smooth $\widetilde f$, $\overline \partial_\rho\kappa={\cal O}
(\rho^\infty)$. 
Now we take real part of (4.27)~:
$$\re p\circ\kappa(x, \xi)=q'_0(\iota')+r'(x, \xi)\leqno(4.28)$$ 
where $\iota'$ stand for the real and imaginary part of 
$\iota$ (it is easy to see that these $2n$ new action variables Poisson
commute for $\{\bullet, \bullet\}_R$.~) 
Following the proof of Theorem 0.1, we consider
the family $q'_s=q'_0+sr'$, $0\leq s\leq 1$. 

As above we look for a family of smooth $\kappa_s$ preserving 
$\re \sigma_{\bf C}$,
satisfying the identity $q'_s\circ \kappa_s=q'_0$, and 
$$\partial_s\kappa_s
=X_s\circ \kappa_s \leqno(4.29)$$
We look for $X_s$ of the form $X_s=\widehat H_{f_s}$,
for some family of real valued functions
$f_s\in I^\infty$. 
Since $q'_s$ is real, we get 
$$\langle \widehat H_{f_s}, q'_s\rangle=-\langle \widehat H_{q'_s}, 
f_s\rangle=-r' \leqno(4.30)$$
and again we are led to solve the 
homological equation $\langle \widehat H_{q'_s}, f_s
\rangle=r'$, for which Proposition 4.3 gives $f_s\in I^\infty$.
Then (4.29) has a solution of the form $\kappa_s=I+\kappa'_s$,
$\kappa'_s\in I^\infty$, uniformly for $s$ on compact sets.
Furthermore, by construction, $\kappa_s$ preserves $\re \sigma_{\bf C}$,
and $(\kappa_s)^*\sigma_{\bf C}=\sigma_{\bf C}$ mod $I^\infty$.
Theorem 0.4 easily follows. $\clubsuit$
\medskip
\noindent {\bf d) Remark : Monodromy along IR-manifolds}.
\medskip
Let $p$ analytic be analytic and have a non degenerate critical point 
at $\rho=0$, such that 
$F_{\rho_0}$ has no purely imaginary, and rationally independent  eigenvalues 
as above. Assume $p$ is real on the real domain. 
We can apply Theorem 0.1 to $T^*{\bf R}^n$ so $p$ is integrable in the
$C^\infty$ sense on the real domain, for some real canonical transform
$\kappa=\kappa_0$ that takes $p$ into its Birkhoff normal form
(we are not interested to know whether $\kappa_0$ is related to
this of Theorem 0.4, since $\kappa$ is only uniquely determined on the 
formal level.)
We set $\Lambda_0=T^*{\bf R}^n$ and try to move $\Lambda_0$ around 
$\rho_0$ in the complex domain, so we consider the family of 
IR-manifolds $\Lambda_s=\exp(isH_p)(\Lambda_0)$, $s\in{\bf R}$ 
(recall that a submanifold of $T^*{\bf C}^n$ is called IR if it is
Lagrangian for $\im \sigma_{\bf C}$ and symplectic for $\re \sigma_{\bf C}$.~)
Then again $p$ is
clearly integrable on $\Lambda_s$, in the $C^\infty$ sense, 
i.e. for real times, and the
problem arises on how far we can go. The 1-d case has been settled
in [HeSj2,App.b], where the authors recover the wellknown fact
that $p$ is integrable in the holomorphic sense~; here  $\kappa$
is univalued, so making a reflection on $\rho_0$ gives $\Lambda_\pi=
\Lambda_{2\pi}=\Lambda_0$. 
This is actually the way that the ``exact Birkhoff normal
form'' was obtained. In several variables we cannnot expect
integrability, nor even recovering $\Lambda_s=\Lambda_0$ for some $s$,
since the orbits may never close (see [Ro2] for a more
complete study of monodromy.~) 

\bigskip
\centerline {\bf Appendix.  The Birkhoff transformations.} 
\medskip
We recall here from [KaRo] some formal constructions, using Lie brackets,
borrowed essentially from
[AbMa,p.500]. There are of course many alternative proofs, 
the idea here is just to 
write formal power series in the most convenient way.  
Since the procedure is mere algebra, it works equally in the
holomorphic, real analytic or 
$C^\infty$ category  (but of course analyticity will be lost during
the game, since small denominators make Birkhoff series divergent, at least in
the common sense. ) In particular eigenvalues can be real or complex.
When eigenvalues are complex, and the hamiltonian real and $C^\infty$,
we can recover real asymptotics just by using an appropriate
linear symplectic transformation of coordinates.
As in [IaSj], we discuss the parameter dependent case. In what follows, 
$s\in\neigh (0, {\bf R}^k)$. 

Let $p=p(s)$ depend smoothly on $s$, and have a 
non degenerate critical point of hyperbolic type at $\rho_s$.
If $p(s)$ is complex valued, we assume also that $\overline
\partial_{(z, \zeta)}p$
vanishes of infinite order at $\rho_s$, so that $p(s)$ 
has formal Taylor series in $(z, \zeta)$ at $\rho_s$.
After a linear symplectic
change of coordinates, depending smoothly on $s$, 
we may assume that $\rho_s=\rho_0=0$, and $p(s)$ has 
quadratic part $p_2(z, \zeta, s)=\Sum_{j=1}^n \lambda_j(s) z_j(s) \zeta_j(s)$. 
We assume also that 
$p(0)$ has rationally independent (or non resonant) frequencies 
$(\lambda_1(0),\cdots,\lambda_n(0))=(\lambda_1,\cdots,\lambda_n)$.
Using that the symplectic group is connected, we may further perform
a symplectic, linear change of coordinates, $C^\infty$ in $s$, such that    
$z_j(s), \zeta_j(s)$ become independent of $s$, and
$p_2(z, \zeta, s)=\Sum_{j=1}^n \lambda_j(s) z_j \zeta_j$.
Of course, the $\lambda_j(s)$'s do not in general
verify the non resonance condition
for $s\neq 0$, but we shall investigate up to which accuracy 
Birkhoff series hold in that case.

After reduction
of the quadratic part as above, $p(s)$ now takes the form
$$p(z, \zeta, s) = p_2(z, \zeta, s) + {\cal O}(|z, \zeta|^3)\leqno (A.1)$$ 
We want to construct a map
$f=f(s)$ between neighborhhods 
${\cal V}(0)$ of
$\rho_0=0 \in T^*{\bf R}^n$, such that $(\exp H_{f(s)})^* H_{p(s)}$ 
is resonant. Indeed we have :
\medskip

\noindent {\bf Proposition A.1}: Let $p(s)=p(z, \zeta, s)$ as above,
and $\rho=(z, \zeta)$. Then
there exists a
smooth canonical transforms $\kappa(s): 
{\cal V}(0) \to {\cal V}(0)$ in
$T^*{\bf R}^n$, and a smooth function $q(s)=q(\iota, s)$, $\iota$
as in (0.3),
such that $\kappa(\rho_0,s) = \rho_0=0$, 
$d\kappa(\rho_0, s)=\id $ 
and
$$p(s)\circ\kappa(s)=q(s)+\rho^3{\cal O}\bigl((\rho, s)^\infty\bigr)$$
 
\smallskip

\noindent {\it Proof}: For simplicity, we assume $k=1$, but the 
general case is similar.
We introduce a small ordering parameter $\e $ and 
rescale coordinates
$(y, \eta)$, as 
$(\e y, \e \eta) = (z, \zeta)$ so that $p(z,
\zeta, s) = \e ^2 p_2(y,
\eta, s) + \e ^3 p_3(y, \eta, s) + \cdots$ where $p_j$ is homogeneous of degree
$j$. Working first at the level of formal Taylor series, we want to solve
(formally), denoting $p=p(s)$, $f=f(s)$~:
$$ (\exp t H_f)^* H_p = \Sum_{j \geq 0} {t^j \over j!}[H_f, [H_f,
\cdots, [H_f, H_p] \cdots ]] = H_r
\leqno(A.2)$$ where $r=r(s)$ is resonant, and $t= \e ^2$.
(Here resonant means of course $H_{p_2}r(s)\sim 0$ (in the sense of Taylor
series at $\rho_0$,~) where $p_2$ as in (A.1), or equivalently
$r(s)\sim r(y_1\eta_1, \cdots, y_n\eta_n, s)$.)
We look also for $f(y, \eta, s) = \e f_1(y, \eta, s) + \e ^2 f_2(y, \eta, s) 
+ \cdots$
with
$f_j$ homogeneous of degree $j+2$.  We proceed by induction.
Collecting the
$\e ^3$-terms in (A.2), we want to find $f_1$ such that $H_{p_3}- H_{\{p_2,
f_1\}}$ is resonant, i.e. 
$p_3 - \{p_2, f_1 \}$ is resonant. Writing $p_3(y, \eta, s) = \Sum_{|\alpha +
\beta|=3}p_{\alpha \beta}(s)y^\alpha \eta^\beta$, 
$f_1(y, \eta, s) = \Sum_{|\alpha + \beta|=3}a_{\alpha \beta}(s)y^\alpha
\eta^\beta$ we try to achieve this condition at any order in $s$.
At zeroth order, i.e. for $s=0$, we take  
$a_{\alpha \beta}(0)=-{ p_{\alpha \beta}(0) \over \langle
\lambda,
\alpha - \beta\rangle }$ for $\alpha\neq\beta$ and $a_{\alpha \beta}(0)=0$
otherwise.
At first order in $s$, the condition that
$\partial_s\bigl(p_3 - \{p_2, f_1 \}\bigr)|_{s=0}$ 
is resonant gives 
$$\partial_s a_{\alpha \beta}(0)={\partial_s p_{\alpha \beta}(0)-
\langle\partial_s\lambda(0),\alpha - \beta\rangle a_{\alpha \beta}(0) \over 
\langle\lambda,\alpha - \beta\rangle }$$
when $\alpha\neq\beta$ and say, $\partial_s a_{\alpha \beta}(0)=0$
otherwise. This process extends by induction to any order in $s$
(note that when $s$ is vector valued, we need to check symmetry
for higher derivatives. )

So far we have constructed the formal Taylor series for $a_{\alpha \beta}(s)$
at $s=0$, and found $f_1(s)$ with an uncertainty $\rho^3{\cal O}(s^\infty)$
(in the original variables).  
Next we collect the $\e ^4$-
terms, which gives :
$$ p_4 - H_{p_2}f_2 - H_{p_3}f_1 + {1\over 2}\{f_1, \{f_1,
p_2\}\} =_{\fed } -H_{p_2}f_2 + q_4$$ 
We want to find $f_2=f_2(s)$ such that
$-H_{p_2}f_2 + q_4$ is resonant. Writing
$q_4(y, \eta, s) = \Sum_{|\alpha + \beta|=4}q_{\alpha \beta}(s)
y^\alpha \eta^\beta$, 
$f_2(y, \eta, s) = \Sum_{|\alpha + \beta|=4}a_{\alpha \beta}(s)y^\alpha
\eta^\beta$ we look again for the Taylor series 
$a_{\alpha \beta}(s)=a_{\alpha \beta}(0)+\partial_s a_{\alpha \beta}(0)s
+{1\over 2}\partial^2_s a_{\alpha \beta}(0)s^2+\cdots$. At zeroth order we
may take $a_{\alpha \beta}(0)= 0$ for $\alpha = \beta$, and
$a_{\alpha
\beta}(0)={ q_{\alpha\beta}(0)\over\langle \lambda, \alpha - \beta\rangle }$
otherwise, then carry on the procedure as above at any order in $s$.
This gives $ H_{p_2}f_2 =r_4$, where $r_4(s) = \Sum_{|\alpha| = 2}
q_{\alpha\alpha}(s)y^\alpha\eta^\alpha$ is the resonant part of $q_4$.  

Assume by induction that we have already constructed $f_1, \cdots, f_{N-1}$
homogeneous of degree $3, \cdots, N+1$, so that $f^{(N-1)} = \Sum_{j=1}^{N-1}
\e ^j f_j$ verifies (A.2) at order $\e ^{N+1}$ in $\rho$, and infinite order 
in $s$. Then we try
$f^{(N)} = f^{(N-1)} + \e ^N f_N$ to fulfill (A.2) up to order $\e ^{N+2}$,
i.e. find $f_N=f_N(s)$ such that 
$$ H_p +  [H_f, H_p]+ \cdots +{t^{N-1} \over (N-1)!}[H_f, [H_f,
\cdots, [H_f, H_p] \cdots ]] + {t^N \over N!}[H_f, [H_f,
\cdots, [H_f, H_p] \cdots ]]$$ is resonant. Each of the terms of
that sum are expanded to order $\e ^{N+2}$. The last one is a $N$-fold bracket
and contains only 
$[H_{f_1}, [H_{f_1},
\cdots, [H_{f_1}, H_{p_2}] \cdots ]]$ to this order~; other terms are $j$-fold
brackets containing $f_1, \cdots, f_{N-1}$, and $f_N$ occurs only in
$[H_{f_N}, H_{p_2}]$. Writing  
$f_N(y, \eta, s) = \Sum_{|\alpha + \beta|=N+2}a_{\alpha \beta}(s)y^\alpha
\eta^\beta$, we can find $a_{\alpha\beta}(s)$ as before so that $-H_{p_2}f_N
+q_{N+2}= r_{N+2}$ where $q_{N+2}=q_{N+2}(s)$ and $r_{N+2}=r_{N+2}(s)$ 
are of degree $N+2$ (or
$r_{N+2}(s)=0$, according to the parity of
$N$), and $r_{N+2}(s)$ is resonant. 

Summing up, we have found $f_j, \ r_j$, deg($f_j) = j+2$, deg($r_j) = j$ such
that 
$$\bigl( \exp tH_{(\e f_1+\e ^2f_2+ \cdots)}\bigr)^*  
H_{(\e ^2p_2+ \e ^3p_3+ \e
^4 p_4 +\cdots)} = H_{\e ^4r_4+ \cdots}$$ so (1.15) is verified at the level of
formal power series. In the original variables
$(z, \zeta) = \e (y, \eta)$, so by  homogeneity~:
$(\exp H_{f(s)})^*H_{p(s)} = H_{r(s)}$. 

All this computation can be implemented at the level of $C^\infty$ germs of
functions at $\rho=(0, 0), s=0$
if we apply Borel's theorem to the $f_j(s)$ and $r_j(s)$.
Hence the relation $(\exp H_{f(s)})^*H_{p(s)} = H_{r(s)}$ 
holds at the level of $C^\infty$
germs, with $r(s)$ resonant, i.e. asymptotic to a $C^\infty$
function of $(z_1\zeta_1, \cdots, z_n \zeta_n)$. Since $(\exp H_f)^*H_p = 
H_{p\circ\exp H_f}$ [AbMa,p.194], we get
$H_{p \circ \exp H_f} = H_r$, and so $p \circ \exp H_f = r$ is
resonant. So we proved the Proposition with $\kappa(s) = 
\exp H_{\widetilde f(s)}$, where $\widetilde f(s)$ is a Borel sum
for $f(z, \zeta, s)$.  
$\clubsuit$

\medskip
\centerline {\bf References}

\bigskip

\noindent [AbM] R. Abraham, J. Marsden. The Foundations of Mechanics. Benjamin,
N.Y. Revised edition, 1985.

\noindent [AbRo] R. Abraham, J. Robbin (with an Appendix of A. Kelley).
Transversal mappings and flows. Cummings, 1967.

\noindent [Ar]  V. Arnold.  M\'ethodes Math\'ematiques de la 
M\'ecanique
classique. Mir, Moscou, 1976. 

\noindent [ArNo] V. Arnold, S. Novikov, {\it eds}. Dynamical Systems III-IV.
Encyclopoedia of Math. Springer, 1988-1990.

\noindent [ArVaGo] V. Arnold, A. Varchenko, S. Goussein-Zad\'e. Singularit\'es
des applications diff\'e- rentiables I. Mir, Moscou, 1986.

\noindent [Au] M. Audin. Les syst\`emes hamiltoniens et leur int\'egrabilit\'e.
Soc. Math. de France, Cours sp\'ecialis\'es {\bf 8}, 2001.
 
\noindent [BaLLWa] A. Banyaga, R. de La Llave, C. Wayne. Cohomological
Equations near hyperbolic points and geometric versions of Sternberg
linearization theorem. J. of Geom. Anal. 6, p.613-649, 1996.

\noindent [Bi] G. D. Birkhoff. Dynamical Systems. Amer. Math. Soc. Colloquium
Publ. 1927, revised ed. 1966.
 
\noindent [Br] F. Bruhat. Travaux de Sternberg. S\'eminaire Bourbaki, no 217,
1960-1

\noindent [Ch] K.T. Chen. Equivalence and decomposition of vector fields...
American J. of Math. 85, p.693-722, 1963. Reprinted in: Collected Papers,
Birkh\"auser, 2001.

\noindent [CuB] R.Cushman, L.Bates. Global aspects of classical 
integrable systems. Birkh\"auser, Basel, 1997.

\noindent [Ec] J. Ecalle. Introduction aux fonctions analysables et
conjecture de Dulac. Hermann, Paris, 1992.

\noindent[E] L. H. Eliasson. Normal forms for hamiltonian systems with Poisson
commuting integrals-elliptic case. Comment. Math. Helv. 65, p.4-35, 1990.

\noindent [Fr] J.P. Fran\c coise. Propri\'et\'es de g\'en\'ericit\'e des
transformations canoniques, {\it in} Geometric dynamics. Proceedings, Rio de
Janeiro. Palis, {\it ed.} 
Lect. Notes in Math. no 1007, Springer, 1983.  

\noindent [Gal] G. Gallavotti. The Elements of Mechanics. Springer, 1983.

\noindent [GiDeFoGaSi] A. Giorgilli, A. Delsham, E. Fontich, L. Galgani, C.
Sim\`o. Effective stability for a Hamiltonian system near an equilibrium 
point... J. of Diff. Eq. 77, p.167-198, 1989.

\noindent [Gr] S. Graff. On the conservation of hyperbolic tori
for hamiltonian systems. J. Diff. Equations 15, p.1-69, 1974

\noindent [GuSc] V. Guillemin, D. Schaeffer. On a certain class of
fuchsian partial differential equations. Duke Math. J. 44(1), p.157-199, 1977.

\noindent [Ha] Ph. Hartmann. Ordinary Differential Equations. Wiley, 1964.

\noindent [HeSj] B. Helffer, J. Sj\"ostrand. {\bf 1.} Multiple wells
in the semi-classical limit III -interaction through non-resonant wells.
Math. Nachr. 124, 1985. {\bf 2.}
Semi-classical analysis
for Harper's equation III. Bull. Soc. Math. de France, M\'emoire
39, Tome 117, 1989.
 
\noindent [IaSj] A. Iantchenko, J. Sj\"ostrand. Birkhoff normal
forms for Fourier integral operators II. Preprint, 2001. 

\noindent [It] H. Ito. Integrable symplectic maps and their 
Birkhoff normal form. T\^ohoku Math. J. 49, p.73-114, 1997.

\noindent [KaRo] N. Kaidi, M. Rouleux. Quasi-invariant tori and
semi-excited states for Schr\"odinger operators I. Asymptotics. 
Comm. Part. Diff. Eq., to appear.

\noindent [LasSj] B. Lascar, J.Sj\"ostrand. Equation de Schr\"odinger \dots
Comm. Part. Diff. Eq. 10(5), p.467-523, 1985.

\noindent [M] P. Malliavin. G\'eom\'etrie diff\'erentielle intrins\`eque.
Hermann. Paris, 1972.
 
\noindent [MaSo] A. Martinez, V. Sordoni. Microlocal WKB Expansions. 
J. of Funct. Anal. 168, p.380-402, 1999.

\noindent [MeSj] A. Melin, J. Sj\"ostrand. Determinants of
pseudo-differential operators and complex deformations of phase space. Preprint
Ecole Polytechnique, 2001.

\noindent [Mo] J. Moser. On the generalization of a theorem of A. Lyapunov.
Comm. Pure Appli. Math. 11, p.257-271, 1958.

\noindent [Ne] E. Nelson. Topics in dynamics I: flows. Princeton University
Press, 1969.

\noindent [Ro] M. Rouleux. {\bf 1.} Quasi-invariant tori and semi-excited 
states for
Schr\"odinger operators II. Tunneling. In preparation. 
{\bf 2.} Integrability of an holomorphic hamiltonian near a hyperbolic
fixed point. In preparation.

\noindent [Si] C. L. Siegel. {\bf 1.} \"Uber die Normalform analytischer
Differential-  gleichungen in der N\"ahe einer Gleichgewichtsl\"osung. Nachr.
Akad. Wiss. G\"ottingen, Math. Phys. Kl., p.21-30, 1952. 
{\bf 2.} \"Uber die Existenz einer Normalform \dots Math. Annalen 128,
p.144-170, 1954.
 
\noindent [SiMo]  C. L. Siegel, J. Moser. Lectures on Celestial
Mechanics. Springer, 1971.

\noindent [Sj] J. Sj\"ostrand. {\bf 1.} Singularit\'es analytiques 
microlocales.
Ast\'erisque 95, Soc. Math. France, 1982. {\bf 2.} Analytic wavefront sets and 
operators with multiple characteristics. Hokkaido Math. J. XII, 3,
p.392-433, 1983. {\bf 3.} 
Functions spaces associated to global I-Lagrangian manifolds, {\it in}
Structure of solutions of differential equations, Katata/Kyoto, 1995,
World Scientific, 1996. {\bf 4.} Semi-excited states in
non-degenerate potential wells. Asymptotic Analysis, 6, p.29-43, 1992.

\noindent [St] S. Sternberg. The structure of local diffeomorphisms III. 
American
J. of Math. 81, p.578-604, 1959

\noindent [Vi] M. Vittot. Birkhoff expansions in Hamiltonian mechanics :
a simplification of the combinatorics, {\it in}~: ``Non-linear
dynamics'' (Bologna, 1988). G. Turchetti, {\it eds.} World Scient.
Singapore, (1989), p.276-286.

\noindent [Vu] S. Vu Ngoc. On semi-global invariants for focus-focus
singularities. Preprint Institut Fourier, 2001

\end